\documentclass{article}

\usepackage{arxiv}

\usepackage[utf8]{inputenc} 
\usepackage[T1]{fontenc}    
\usepackage{hyperref}       
\usepackage{url}            
\usepackage{booktabs}       
\usepackage{amsfonts}       
\usepackage{nicefrac}       
\usepackage{microtype}      
\usepackage{lipsum}		
\usepackage{graphicx}
\usepackage{natbib}
\usepackage{doi}

\usepackage{amssymb}
\usepackage{amsmath}
\usepackage{amsthm}

\newtheorem{problem}{Problem}
\newtheorem{theorem}{Theorem}
\newcommand{\ceil}[1]{\left\lceil #1 \right\rceil}

\usepackage{forest}
\usepackage{multirow}
\usepackage{longtable}
\usepackage{booktabs}
\usepackage[tableposition=below]{caption}
\usepackage{comment}

\LTcapwidth=\textwidth  
\title{A Hierarchical Integer Linear Programming Approach for Optimizing Team Formation in Education}

\author{ 
\href{https://orcid.org/0009-0008-2013-1709}{\includegraphics[scale=0.06]{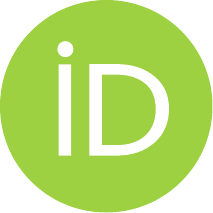}\hspace{1mm}Aaron Kessler}
\\
    Department of Computer Science \\
    Trier University of Applied Sciences\\
	Trier, 54293 \\
    Rhineland-Palatinate, Germany\\
	\texttt{kesslea@hochschule-trier.de} \\
	\And
	\href{https://orcid.org/0009-0008-3065-5178}{\includegraphics[scale=0.06]{orcid.pdf}\hspace{1mm}Tim Scheiber} \\
    Department of Computer Science \\
    Trier University of Applied Sciences\\
	Trier, 54293 \\
    Rhineland-Palatinate, Germany\\
	\texttt{scheibet@hochschule-trier.de} \\
    \And
	\href{https://orcid.org/0000-0002-3371-6983}{\includegraphics[scale=0.06]{orcid.pdf}\hspace{1mm}Heinz Schmitz} \\
    Department of Computer Science \\
    Trier University of Applied Sciences\\
	Trier, 54293 \\
    Rhineland-Palatinate, Germany\\
	\texttt{h.schmitz@hochschule-trier.de} \\
     \And
	\href{https://orcid.org/0000-0002-4243-4128}{\includegraphics[scale=0.06]{orcid.pdf}\hspace{1mm}Ioanna Lykourentzou} \thanks{Corresponding author}\\
    Department of Information and Computing Sciences\\
    Utrecht University\\
	Princetonplein 5, 3584 CC \\
    Utrecht, The Netherlands\\
	\texttt{i.lykourentzou@uu.nl} \\
}

\renewcommand{\headeright}{}
\renewcommand{\undertitle}{Preprint}
\renewcommand{\shorttitle}{A Hierarchical Integer Linear Programming Approach for Optimizing Team Formation in Education}

\hypersetup{
pdfkeywords={team formation, multi-objective optimization, integer linear programming, problem taxonomy, education},
}

\begin{document}
\maketitle

\begin{abstract}
Teamwork is integral to higher education, fostering students' interpersonal skills, improving learning outcomes, and preparing them for professional collaboration later in their careers.
While team formation has traditionally been managed by humans, either instructors or students, algorithmic approaches have recently emerged to optimize this process.
However, existing algorithmic team formation methods often focus on expert teams, overlook agency in choosing one’s teammates, and are limited to a single team formation setting.
These limitations make them less suitable for education, where no student can be left out, student agency is crucial for motivation, and team formation needs vary across courses and programs.
In this paper, we introduce the EDUCATIONAL TEAM FORMATION problem (EDU-TF), a partitioning optimization problem model tailored to the unique needs of education, integrating both teacher and student requirements.
To solve EDU-TF, we propose a modular optimization approach, one of the first to allow the flexible adjustment of objectives according to educational needs, enhancing the method's applicability across various classroom settings rather than just research environments.
Results from evaluating ten strategies derived from our model on real-world university datasets indicate that our approach outperforms heuristic teacher-assigned teams by better accommodating student preferences.
Our study contributes a new modular approach to partition-based algorithmic team formation and provides valuable insights for future research on team formation in educational settings.
\end{abstract}


\section{Introduction}
Teamwork is an important aspect of many higher education programs.
Collaborating on team projects stimulates cooperative learning, fosters shared goals and mutual learning responsibility among students, enhances interpersonal skills, and deepens the students’ understanding of the material \citep{ravenscroft1995incentives, smith1992collaborative, van2007impact, johnson2000cooperative}.
Compared to traditional teaching methods, students engaged in collaborative learning tend to retain information longer, exhibit lower dropout rates, and gain a better understanding of real-world work environments \citep{oakley2004turning}.
Performance-wise, team projects consistently improve students' academic achievements over more conventional instructional approaches \citep{oakley2004turning, springer1999effects, ahmad2010effects}.

The team with which one collaborates plays a critical role in achieving positive outcomes in collaborative learning \citep{cohen1997makes}.
An effective team optimally leverages individual skills while complementing them with those of other members, enhancing both learning and performance.
Additionally, working with people who share good rapport increases psychological safety, reduces unproductive conflicts, and improves overall learning outcomes, satisfaction, and team viability.
However, forming optimal teams -- considering both skill alignment and personal compatibility -- remains a significant challenge due to the instructor’s limited time, lack of sufficient knowledge about individual students' preferences, and capacity constraints, especially in medium- and large-sized university courses.

Existing heuristic team formation methods employed by university instructors often fall short.
\textit{Instructor-led team formation} \citep{hilton2010instructor, rusticus2019comparing}, where the teacher assigns students to teams based on criteria such as skill alignment or project preferences may ensure diverse enough teams to meet course objectives, however, it can lead to interpersonal conflicts, requests for team changes, suboptimal learning outcomes, demotivation, and even dropouts \citep{akbar2018improving}.
\textit{Student-led team formation}, where students self-select their teammates, is another common approach.
Teams formed in this manner are typically based on friendships or past collaborations, and members often share similar backgrounds, interests, and commitment levels \citep{hilton2010instructor}.
While this method fosters autonomy, strengthens relationships, and can even improve performance \citep{chapman2006can}, it has significant drawbacks.
Self-formed teams may lack the necessary skills and expertise for the task at hand \citep{basta2011role, chapman2006can}.
Additionally, this approach tends to lead to the formation of ``all-star'' teams, where stronger or more popular students group together, leaving behind those who are weaker, from underrepresented backgrounds, or less integrated into the social network.
This segregation results in unequal learning opportunities and subpar performance \citep{oakley2004turning, aranzabal2022team, bacon2001methods, miglietti2002using}.
The third method, \textit{random team assignment}, involves a purely arbitrary process.
Although easy to implement, particularly in large courses such as MOOCs, this method fails to optimize team formation, leading to significant variations in team performance and perceptions of unequal effort \citep{pociask2017does}.

Algorithmic team formation has emerged as a solution to address the limitations of manual human approaches, typically employing exact techniques such as constraint programming.
However, existing team formation algorithms usually exhibit three shortcomings that make them unsuitable for educational settings.
First, most team formation algorithms focus on expert team selection, where the objective is to form selected teams from a pool of experts \citep{bahargam2019team, bhowmik2014submodularity}, without the requirement that every participant is assigned to a team.
These algorithms, oftentimes designed for business environments to facilitate recruitment and batch-hiring \citep{golshan2014profit} and optimize the allocation of human resources to projects \citep{gutierrez2016multiple, ballesteros2019human}, are unsuitable for educational settings, where all students must be included in the team formation solution.
Second, the large majority of team formation algorithms operate in a top-down manner, taking into account only the objectives of the manager or teacher (e.g., regarding skill coverage and team size), while ignoring participant preferences.
However, participant agency, most notably one’s ability to have a say in who they work with, is critical for motivation and overall teamwork success, especially in educational settings.
Third, most of the currently proposed algorithms are designed to optimize a highly specific setting, making them impractical for use in real-world educational environments, where team formation needs vary per program and course, an issue further complicated by the nascent and non-systematic nature of the literature, which often focuses on specific problem definitions rather than a generic approach, making it difficult to compare, generalize, and apply the proposed algorithms.

In this paper we make the following contributions:
\begin{enumerate}
    \item We take a ``bird’s eye'', yet comprehensive, scoping view of the existing algorithmic team formation literature and propose a systematic taxonomy of the different problem families within it, making it easier to locate the specific problem family being worked on.
    The proposed taxonomy is defined by the number of teams that a candidate can join and the number of tasks considered in the problem setting.
    We distinguish between the Top-k team formation and Partition team formation families in this taxonomy, applied to both single and multiple tasks.
    We differentiate between problem families (i.e., generic problem settings) and problem definitions that capture specific tasks with particular attributes, inputs, and constraints. We demonstrate that the EDUCATIONAL TEAM FORMATION (EDU-TF) problem, which is the focus of this paper, is a problem definition derived from the Single-task Partition TFP family.
    \item We formalize the EDU-TF problem, defined specifically for the requirements of educational contexts.
    It is among the first to incorporate not only top-down teacher constraints (for skill coverage and team size), but also bottom-up student objectives (for preferred teammates) into the optimization problem definition.
    We analyze the computational complexity of the EDU-TF problem and show it is \textbf{NP}-hard, a result derived from the fact that its decision variant that ignores teammate preferences and only asks for the existence of a feasible solution is already \textbf{NP}-complete.
    \item We introduce a hierarchical integer linear programming (ILP) approach for solving the EDU-TF problem.
    Our approach is among the first to allow teachers to flexibly adjust and prioritize the optimization objectives of team formation according to their needs, making it practically applicable to a variety of classroom settings.
    \item We experimentally validate our approach by comparing ten team formation strategies derived from it with the heuristic teacher solution across nine real-world datasets.
    Our results show that our approach outperforms the heuristic teacher solution in accommodating student preferences while being sufficiently computationally efficient, achieving both optimal solutions and good-enough solutions within 15 minutes.
\end{enumerate}

The rest of this paper is organized as follows.
Section \ref{sec:taxonomy} reviews related work, categorizing it into a structured taxonomy of team formation problem families, commencing from the General Team Formation Problem.
Section \ref{sec:formalization} formalizes the EDU-TF problem and analyzes its complexity.
Section \ref{sec:ilp} presents our modular ILP to solve EDU-TF, which relies on a base model with the objectives as extensions that can be combined and (re-)ordered according to the specific teacher needs.
Section \ref{sec:results} presents experimental results, comparing ten different strategies derived from our modular approach with the respective teacher heuristic solution on nine different datasets derived from eight real university courses.
Section \ref{sec:discussion} discusses the contributions, limitations, and future work, and Section \ref{sec:conclusion} concludes the paper.

\section{The Team Formation Problem (TFP): Taxonomy of problem families}
\label{sec:taxonomy}

\begin{figure}
    \includegraphics[width=\linewidth]{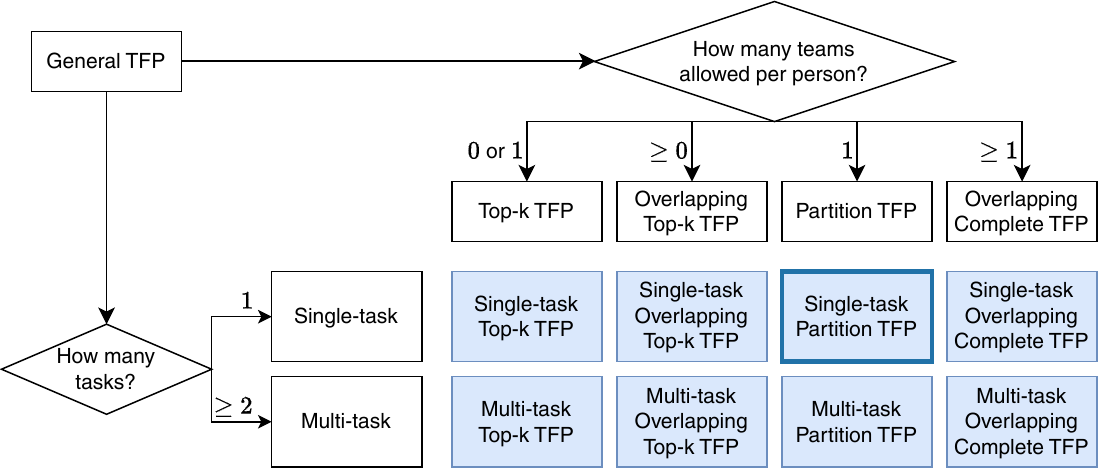}
    \caption{
        Taxonomy of the General Team Formation Problem (General TFP) families, defined by the number of teams a candidate can join and the number of tasks considered in the problem setting.
        The Top-k TFP, the earliest variant explored in TFP literature, models the selection of one or more expert teams.
        This problem family allows for excluding candidates and is commonly applied in business settings.
        The Partition TFP divides the pool of candidate team members into teams, ensuring each person is allocated to exactly one team.
        Overlapping problem families extend both Top-k and Partition TFPs by allowing individuals to join multiple teams.
        From these problem families, various problem definitions can be derived, concretizing the problem family according to specific attributes and modeling elements (inputs, constraints, etc.).
        EDUCATIONAL TEAM FORMATION (EDU-TF), on which we focus in this paper, is a problem definition derived from the Single-task Partition TFP family, in which no student can be left out, each student must belong to exactly one team, and all teams work on the same task.
    }
    \label{fig:problem_taxonomy}
\end{figure}

The General Team Formation Problem (General TFP) can be understood as the challenge of identifying one or multiple teams of collaborators to solve one or multiple tasks from a pool of candidate teammates, in such a way as to fulfill specific requirements (individual, team, or organizational ones) as best as possible.
The General TFP comprises several distinct problem families, depending on the number of teams a person can be part of and the number of tasks these teams will handle (see Figure~\ref{fig:problem_taxonomy}).

Starting with the number of teams to be formed, we find four main families.
The \textit{Top-k TFP} variant, by far the one mostly explored in the literature, pertains to selecting a number of people -- usually referred to as experts -- from a pool of candidate teammates and assigning them to work, in one or more teams, on a single task (Single-task Top-k TFP) or multiple ones (Multi-task Top-k TFP).
The Top-k TFP is characterized by the possibility of excluding some candidates from the team formation solution, with each selected team member belonging to one team at most.
This variant is commonly studied by researchers through the case of expert team selection in social networks and applied in practice in business settings, where workers are allocated to client projects but can also remain ``on the bench''.

The second most explored variant, though still far less studied than the Top-k TFP, is the \textit{Partition TFP}.
Here, the goal is to allocate every candidate team member to exactly one team, thus partitioning the candidate pool into disjoint subsets.
In this variant, no person can be left out of the team formation solution and each person belongs to exactly one team.
This variant is frequently applied in educational settings, where students must be assigned to project-based course teams.
Depending on whether the teams work on the same task or on multiple distinct tasks, this problem variant can be further divided into the Single-task Partition TFP or the Multi-task Partition TFP.

The next two problem families are included in our taxonomy by means of inference, as they are less explored in team formation optimization literature, despite their frequent use in practical settings and occasional discussion in social science research.
These are the \textit{Overlapping Top-k TFP} and the \textit{Overlapping Complete TFP}.
Overlapping means that a person, if selected, can belong to multiple teams.
In the Overlapping Top-k TFP, similar to the traditional Top-k TFP, a person can be part of multiple teams or none at all.
In contrast, the Overlapping Complete TFP, akin to the Partition TFP, requires that every individual be part of at least one team but allows them to participate in multiple teams simultaneously.
Both problem families are commonly applied in exploratory team formation contexts, such as hackathons or early-stage innovation contests, where participants may engage with multiple teams until project completion or until they choose to leave a team.

In the following, we present the literature supporting this taxonomy and position the EDUCATIONAL TEAM FORMATION Problem within it.

\subsection{Top-k TFP}

The \textit{Top-k TFP} describes the problem where $k$ teams need to be selected out of a pool of candidate participants to complete one or more given tasks, not all participants need to pertain to a team (i.e., participants can be excluded), and there is no person-to-team overlap (i.e., one person can only belong to one team). 

The \textit{Single-task Top-k TFP} refers to the problem where one or more team(s) of experts are selected to work on \textit{a single task}.
Among the first works exploring this problem is the study of \citet{lappas2009finding}, who model the problem of identifying the best team of experts ($k$=1) to solve a given task, aiming to ensure that the team possesses all required skills for the task and has minimal communication costs over a social network.
Generalizing this modeling, \citet{5590958} associated each required skill with a specific number of experts, and proposed the Enhanced-Steiner algorithm to solve this problem instance.
Recognizing that skill is only one of the criteria required to build an excellent team, \citet{zhang2013multi} proposed a multi-objective optimization model that also considers interpersonal relations, and proposed an algorithm offering alternative solutions on a Pareto front to solve it.
Additional examples of studies that focus on the Single-task Top-k TFP are those by \citet{bhowmik2014submodularity}, who incorporate a multitude of team formation criteria using an unconstrained submodular function maximization problem modeling, as well as those by \citet{majumder2012capacitated}, \citet{durfee2014using}, and \citet{farhadi2011effective}, who propose more generalized versions of the problem and test it with larger datasets.

Another early study exploring the Single-task Top-k TFP problem variant is the study of \citet{kargar2011discovering}, who work on the problem of identifying multiple teams of experts to work on a given task, while ensuring that each team has the required skills and the intra-team communication costs are minimized over a social network.
The study proposes a procedure that produces the Top-k teams with a polynomial delay.
Multiple researchers have followed up on this work, accounting for additional attributes besides skill, different cost functions, and additional team collaboration scenarios \citep{rahman2015task, kou2020efficient, gajewar2012multi, ballesteros2019human}.
However, up until 2021, the problem of forming $k$ expert teams working on a single task in polynomial time was still an open problem when the number of expert skills (or, more generally, the individual person attributes) that needed to be accounted for was greater than two.
In 2021, \citet{lokshtanov2021eth} proposed an Exponential Time Hypothesis (ETH)-tight algorithm to form multiple teams of experts and account for multiple individual attributes.
The algorithm runs in polynomial time up to $d = c\log\log(n)$ many attributes per candidate where $n$ is the number of candidate teammates and $c>0$ is a sufficiently small constant.

Departing from the problem variant where the selected team(s) all work on a single task, research has also explored the \textit{Multi-task Top-k TFP} variant, where the team(s) need not only be selected but also to be matched to different tasks.
In this context, the work of \citet{golshan2014profit} studies the problem of hiring teams of experts that collectively possess the skills to solve as successfully as possible as many tasks as possible from a company task collection, without exceeding a given budget. 

In the years that followed, multiple researchers studied extensively the Top-k TFP, most notably from the database, data mining and web research communities, exploring linear programming \citep{campelo2020sociotechnical, rostami2023deep}, goal programming \citep{dacs2022set}, stochastic programming \citep{rahmanniyay2019multi}, matrix factorization \citep{wu2021task}, fuzzy optimization \citep{kalayathankal2021modified}, game theory \citep{sorkhi2012game}, facility location analysis \citep{neshati2014expert, neshati2012multi}, as well as heuristic and meta-heuristic approaches \citep{jiang2019group, gutierrez2016multiple}, to name just a few.
\citet{wang2015comparative} and \citet{juarez2021comprehensive} present two comprehensive state-of-the-art reviews of the different problems and algorithmic approaches to solve them.
Overall, the majority of research works focusing on the Top-k TFP mainly apply their findings to the domains of social networks, crowdsourcing, and business.

\subsection{Partition TFP}
The Partition TFP consists of distributing a set of individuals into disjoint teams so that no individual is left without a team and each person belongs to exactly one team.

Whereas the Top-k TFP is commonly encountered in business settings to facilitate hiring, the Partition TFP often finds application in education, where teachers must divide a project-based course class into teams, ensuring that no student is left out.
Educators typically solve this problem manually, relying on heuristics based on their experience \citep{dzvonyar2018team, lingard2002teaching, bosnic2013picking, sedano2016green, lago2009designing}.
Although this variant has been explored far less than the Top-k TFP, research on it has recently expanded to both single- and multi-task scenarios.
EDUCATIONAL TEAM FORMATION (EDU-TF) belongs to the Partition TFP category.
In this paper we specifically focus in its single-task version.

The first variant is the \textit{Single-task Partition TFP}, where all teams work on the same task.
\citet{andrejczuk2019synergistic} model the scenario where students must be partitioned into teams that possess the necessary competencies and exhibit high diversity in terms of gender and personality.
This study proposes a linear programming algorithm for smaller problem instances and a heuristic for larger ones.
\citet{agrawal2014grouping} model the scenario where students need to be partitioned into non-overlapping groups using two objective functions: one that maximizes team ability (favoring stronger students) and another that maximized the number of less-skilled students placed in teams with more skilled peers.
They propose an efficient approximation algorithm, the complexity of which grows logarithmically with respect to team size.
More recently, 
\citet{candel2023integer} model a scenario where students are assigned to teams to maximize a function quantifying team performance through personality matching.
Their model also incorporates constraints regarding which student pairs should or should not be placed together and proposes a linear integer programming model to solve it.
This study is among the few that consider student-preferred and non-preferred teammates, though it does so using binary preferences defined by the teacher rather than a range of preference values specified by students, as in our model.
Beyond the educational domain, \citet{bahargam2019team} study the problem of partitioning teams to minimize faultline scores, i.e., hypothetical divisions that split a participant pool into relatively homogeneous subteams based on multiple attributes, and propose an efficient heuristic algorithm to solve it.

The second variant is the \textit{Multi-task Partition TFP}, where teams work on different tasks.
Even fewer studies focus on this variant.
\citet{dzvonyar2018algorithmically} examine a scenario in which students in a software project course are divided into teams, each team working on a different project.
The single optimization objective is to assign students to a project as high as possible on their preference list.
Comparing algorithmic and manual team formation, this study shows that their proposed linear programming model and algorithm improved the mean project priority for more than half of the final teams.
\citet{georgara2022allocating} generalize the problem studied by \citet{andrejczuk2019synergistic} by modeling the allocation of multiple agents (e.g., students) to multiple tasks, ensuring no agent-to-team or team-to-task overlaps.
They show that the problem is \textbf{NP}-complete and propose a heuristic solution. 
However, their model does not explicitly enforce the constraint that all agents must belong to a team, even though this requirement is implied in the problem description.
Finally, \citet{chen2025multi} tackles multi-task team formation in an online collaborative learning setting, in which participants collaborate on different learning tasks through a central online platform.
The authors highlight the importance of synchronous online communication and the benefits of collaboration for the development of the participants' skills.
Considering these components, they partition participants and assign them to tasks using a Monte Carlo Tree Search algorithm.

\subsection{Overlapping TFP}
By inference, we can distinguish two additional problem families: \textit{Overlapping Top-k TFP} and \textit{Overlapping Complete TFP}.
In these, each person can, if selected, participate in multiple teams simultaneously and, depending on whether it is the single- or multi-task version, work on multiple tasks as well.
This contrasts with the exclusive task and/or team allocation found in the Top-k and Partition TFP families.
Overlapping TFP is not typically applied in education but finds applications in hackathons or other bottom-up collaboration settings, where individuals contribute to multiple projects and teams to maximize their profit.

To our best knowledge, the \textit{Single-task Overlapping Top-k TFP}, or the \textit{Overlapping Complete TFP} have not yet been examined in the optimization literature.
Three studies examining the \textit{Multi-task Overlapping Top-k TFP} are the ones by \citet{gutierrez2016multiple}, \citet{campelo2021integer}, and \citet{selvarajah2021unified}.
The first study models the problem of allocating multiple people, either full-time or part-time, on multiple project teams based on the skill requirements of each project.
They propose three algorithms to solve the problem, based on Variable Neighbourhood Search, Constraint Programming and Local Search, with the Variable Neighborhood Search algorithm slightly outperforming the other two.
The second study models a similar setting and proposes an Integer Linear Programming (ILP) approach to solve it.
The third work proposes a general framework for selecting experts in a social network to work on a set of projects, which considers the candidates' expertise, communications costs, and geographical proximity among other factors.
While the framework does not explicitly foresee that an expert may be assigned to multiple teams, this scenario is depicted in their examples where the algorithm allows individual experts to be allocated to multiple teams.
They propose a multi-objective evolutionary algorithm to optimize the project teams, which, given the competing nature of certain objectives, returns a Pareto front of possible team set solutions.

\subsection{TFP modeling attributes and the importance of agency}

Research has examined team formation modeling using various attributes, often specific to the application domain.
While these attributes can essentially be seen as either node or edge values in a graph model of potential teammates, they are crucial for two reasons.
First, the attributes help differentiate the problem families presented above, enabling the development of concrete problem definitions (input data, constraints, output).
Second, from a user and application perspective, understanding which attributes are used in the problem definition enhances the explainability of the constructed definitions.
This is particularly important in education and decision-making, as it facilitates better evaluation of the problem definition's effectiveness.
Below, we provide an overview of the most commonly used attributes in the team formation optimization literature, with a specific focus on the role of agency, expressed as direct teammate preference input from the candidates themselves, which team formation algorithms must account for.

\textit{Skill}, in terms of individual knowledge, expertise, roles, communication competencies, or leadership competencies, is one of the primary attributes of team formation explored in the literature.
Skill-based optimization is used in workplace settings to inform worker selection and recruitment \citep{lele2015formation}, to hire crowd worker teams on online labour market platforms \citep{anagnostopoulos2018algorithms, hamrouni2020optimal}, and to build agile teams \citep{rostami2023deep} for software development.

The objective of the team formation algorithm is usually to ensure that the selected teams cover the skills necessary to achieve the team goal(s) \citep{yannibelli2012deterministic}, 
but other objectives such as inter-team skill homogeneousness and intra-team skill heterogeneousness may also be envisioned \citep{moreno2012genetic}. 
Researchers have also focused on ensuring smooth \textit{communication} among the team members.
Communication, also referred to as team member compatibility, interaction quality, or social cohesion, has long been identified in sociological research as a core element of teamwork success \citep{moreno1941foundations}. Previous small group research suggests that cohesive groups tend to perform approximately 18\% better than non-cohesive groups \citep{evans1991group}. In this line \citep{anagnostopoulos2012online, kargar2011discovering, lappas2009finding, rangapuram2013towards, berktacs2021branch} all focus on ensuring that the team of experts selected to work on a specific task has a low communication overhead, in addition to possessing the necessary skills and dividing the workload fairly. 

\textit{Personality} is another frequently used attribute, as prior research indicates that certain personality combinations may perform more efficiently together, due to fewer frictions or power struggles \citep{lykourentzou2016personality}.
In this line \citet{sanchez2023comparing} and \citet{candel2023integer} used personality, as an individual user (i.e., node) attribute, to form student teams using the Belbin and Myer-Briggs personality tests, and \citet{andrejczuk2018solving} used Post-Jungian Personality Theory, respectively.

Adapting team \textit{diversity} is a fourth factor researchers examine as a group formation measure. Prior research suggests that average group heterogeneity, i.e., neither extremely homogeneous nor extremely diverse groups, can increase group performance by up to 10\% \citep{gibson2003healthy, chen2017too}. Under this light, \citet{miranda2020multi} maximize the number of students with different profiles per team in an educational setting, and \citet{entani2022group} examine the problem of minimizing or maximizing diversity among the selected team members' self-rated profiles in the context of collaborative learning. 

Other factors that have been considered include the minimization of ~\textit{budget}, expressed as the costs of hiring a team member \citep{golshan2014profit, kargar2011discovering}, taking into account the participants' \textit{time availability} \citep{durfee2014using}, and balancing the \textit{workload} among the selected team members \citep{anagnostopoulos2012online}.

\paragraph{The role of user input in team formation optimization}
As the reader may have noted, most team formation algorithms do not take into account any direct input from the candidate teammates.
However, research shows that fixing teams in a purely top-down manner, considering only the decision maker's (e.g., teacher, manager) requirements, can stifle team members' creativity, initiative-taking, and the teams' ability to adapt their problem-solving strategies \citep{retelny2017no}, especially for creative tasks.
Conversely, as studies in accountable governance show, accounting for \textit{agency}, expressed as having a high degree of autonomy in defining one's work environment including who one works with, improves intrinsic motivation, increases the sense of ownership over one's own results and the results of their teamwork, promotes creativity, improves satisfaction and oftentimes performance \citep{haas2016secrets, seibert2011antecedents, spector1986perceived, lorenzen2009creativity, andriopoulos2001determinants}.
Specifically for the educational domain, students expressed the need to have a strong voice in the criteria used by the team formation algorithm \citep{jahanbakhsh2017you}.
Despite this, very few optimization algorithms proposed so far consider student input as a modeling attribute.
The study by \citet{asee_peer_44448} is a step in this direction, proposing a team formation optimization model that considers not only student skills and teacher requirements but also individual preferences regarding which projects (but not which teammates) they would like to work on.
Even fewer studies consider the students' preferences regarding which teammates they want to work with, despite the fact that the lack of agency over who one works with is widely documented to increase psychological discomfort \citep{rasmussen2006teamwork}, reduce autonomy, and alienate team members \citep{lawler2009designing, de2001minority}.
In this context, \citet{hastings2020lift} propose LIFT, a novel learner-centered workflow where students propose, vote for, and weigh the criteria used as inputs to the team formation algorithm.
Results from this study showed that students valued being taken into account and that the criteria they chose differed from those frequently chosen by instructors, such as grade (GPA), with students preferring criteria promoting efficient collaboration.
Finally, the recent studies by \citet{gomez2024augmenting} and \citet{gomez2024assigning} involving students, as well as academic staff members, showed that agency has a significantly positive impact on team performance not as a standalone attribute but when combined with diversity in terms of demographics and skills 
\citep{gomez2024augmenting, gomez2024assigning}.

In summary, our study differs from and complements existing literature in three ways.
First, most studies on algorithmic team formation focus on expert teams, where individuals can be excluded from the solution, and thus find application in domains such as business and crowd work.
Fewer studies focus on the partitioning problem, which is relevant to education, where no student can be left out of the solution; our study is one of those.
Second, while skill is widely recognized as one of the most important characteristics to consider, student preferences are rarely taken into account, although agency in selecting one's teammates is a well-documented factor that critically affects teamwork, especially in student-centered educational settings.
Our paper is among the first to structurally incorporate this factor into the problem modeling.
Third, the majority of algorithms proposed in the existing literature function within a specific problem setting, making them less suitable for real-world educational settings where teachers have different requirements per program and even per course.
In contrast, our model is one of the first to permit flexible adjustment of objectives by the teacher based on their needs, thus allowing its use in practice rather than in purely research environments and enabling its generalization to multiple classroom settings.

\section{Formalization of educational team formation}
\label{sec:formalization}

We start the formal model for team formation in an educational setting with a decision problem that asks for feasible solutions, where students are partitioned into teams of given size such that the required skill coverage per team is ensured.
After that, we define different objectives that additionally optimize various kinds of interpersonal preferences of students, extending the initial decision problem into a family of optimization problems that capture desired characteristics of educational team formation.

\paragraph{Notation}
We consider $0 \in \mathbb{N}$ and $\mathbb{N}^+ = \mathbb{N} \setminus \{0\}$.
Let $[x,y]=\{ i \in \mathbb{Z} ~|~ x \le i \le y\}$ for $x,y\in \mathbb{Z}$ with $x\le y$, and $[x]=[0, x-1]$ for $x\in \mathbb{N}^+$.
$\#A$ denotes the cardinality of any finite set $A$.
We call a partition $P$ an $n$-partition if $\#P = n$.

\subsection{Input data and feasible solutions}
\label{sec:problem-data}

At its core, team formation in the educational setting is typically concerned with partitioning a group of students into a desired number of teams of a roughly equal number of persons (\emph{team-size}).
We focus on the case where teams are given similar tasks that all require the same set of skills.
Hence each team must comprise a mix of students that provide, in total, the skills required to perform the task effectively (\emph{team-skill}).
However, a broad range of skills within a team does not automatically translate to the best possible outcome. 
Equally important is the chemistry and cohesion among team members.
This can be expressed by \emph{student preferences} that quantify to which extent a student would like 
to work with any other student. 
In more formal terms, the input to our model consists of the following data:

\begin{itemize}
    \item \textbf{Team data.}
    A course consisting of $m \in \mathbb{N}^+$ \emph{students} must be partitioned into $n \in \mathbb{N}^+$ \emph{teams}.
    The additional \emph{size bounds} $k_{min}, k_{max} \in \mathbb{N}^+$ restrict the size of each team.
    We assume w.l.o.g.~$k_{min} \le m / n \le k_{max}$ since otherwise no feasible $n$-partition of $m$ students exists within the size bounds.
    
    \item \textbf{Skill data.}
    Each course is associated with a finite set of \emph{skills} $S \subset \mathbb{N}$ of interest.
    Every student may possess any of these skills, resulting in a personal \emph{skill set} $S_a \subseteq S$ for each student $a \in [m]$.
    We say that student $a$ \emph{covers} skill $i$ if $i \in S_a$.
    The notion of \emph{skill coverage} extends to teams, i.e., a team covers a skill if at least one of its members covers that skill.
    The \emph{minimum skill coverage} $c \in [0, \#S]$ states the minimum number of skills each team has to cover in a feasible solution.
    
    \item \textbf{Preference data.}
    For each student $a$ and possible teammate $b \in [m] \setminus \{a\}$ the preference value 
    $p_{a,b} \in [-d, d]$ for a fixed $d \in \mathbb{N}^+$ reflects $a$'s desire to collaborate in a team with $b$.
    Higher values indicate higher preferences, while $0$ represents a neutral position. 
    Negative values correspond to an active desire by $a$ not to work with $b$.
    The reflexive preferences $p_{a,a}$ are fixed to $0$.
    Preference values are summarized in a \emph{preference matrix} $P \in [-d, d]^{m \times m}$, which is not symmetric in general.
\end{itemize} 

A \emph{feasible solution} for this input data  is an $n$-partition $\mathcal{T}$ of $[m]$ such that for every team $T \in \mathcal{T}$ it holds that
\begin{equation}
    k_{min} \le \#T \le k_{max} \quad (\emph{team-size constraint})
    \label{eq:team-size-constraint}
\end{equation}
and 
\begin{equation}
    \#\left(\bigcup_{a \in T} S_a\right) \ge c \quad (\emph{team-skill constraint}).
    \label{eq:team-skill-constraint}
\end{equation}
Note that for $c = 0$ there always exists a trivial feasible solution respecting (\ref{eq:team-size-constraint}).
This is important for practical reasons in case the given input does not allow for a feasible solution with $c>0$ at all, but the course has to start with some partition anyway.

\subsection{Feasiblity of educational team formation is \textbf{NP}-complete}

The first algorithmic question is a feasibility problem: 
Given the input, does a feasible solution exist?

\begin{problem}
    EDUCATIONAL TEAM FORMATION (EDU-TF)\\
    Given team and skill data, is there an $n$-partition of students such that every team respects the team-size constraint (\ref{eq:team-size-constraint}) and team-skill constraint (\ref{eq:team-skill-constraint})?
\end{problem}

Since the length of a solution is polynomially bounded by the length of the input and the constraints can be verified in polynomial time, EDU-TF is a member of \textbf{NP} (for basic notions  of computational complexity we refer to a standard textbook, e.g., \citet{gol10}).
Next we show completeness for \textbf{NP} by reduction from SET COVER, a well-known \textbf{NP}-complete problem \citep{Karp1972}.
It asks whether we can select a restricted number of given subsets such that their union is equal to the universe.

\begin{problem}
    SET COVER\\
    Given a finite set $U \subset \mathbb{N}$ (the universe), a family $\mathcal{S}$ of subsets of $U$ and some $k \in \mathbb{N}^+$, is there a subset $\mathcal{C} \subseteq \mathcal{S}$ with $\#\mathcal{C} \le k$ and $\bigcup_{S_i \in \mathcal{C}} S_i = U$?
\end{problem}

W.l.o.g. we make some assumptions about the instances of SET COVER considered here:
Cases where $k = 1$ can be solved in polynomial time simply by checking if any of the given subsets is equal to $U$.
The same is true if $\#\mathcal{S} = 1$.
So from now on suppose 
$\#\mathcal{S}\geq 2$ and $2 \leq k \leq \#\mathcal{S}$.

Since we want to reduce a selection problem to a partition problem, where no entity can be excluded from a solution, we will set up the reduction in a way that all but one of the $n$ teams can be formed trivially.
This is achieved by adding $n-1$ so-called \emph{all-rounders} to the instance, i.e., students that cover all skills. 
Then there is a SET COVER solution if and only if also the remaining team can be formed in a way that all skills ($=U$) are covered.
Here is a small example to illustrate this:
Let a SET COVER instance be given with subsets $\mathcal{S} = \{\{1, 2\}, \{3\}, \{3, 4\}\}$ over the universe $U = \{1, 2, 3, 4\}$, and we are only allowed to select at most $k = 2$ of them.
Then we will have three students with respective skill sets $\{1, 2\}, \{3\}, \{3, 4\} \subseteq S = U$, team size bound $k_{max} = k = 2$ and minimum skill coverage $c = 4 = \#U$.
Clearly, one team of the EDU-TF solutions corresponds to the desired SET COVER solution and the other team becomes feasible if we add an additional student with skill set $\{1, 2, 3, 4\}$ to the instance:
\begin{equation*}
     \underbrace{\{ \{1, 2\}, \{3, 4\} \}}_\text{team 1}, \underbrace{\{ \{3\}, \overbrace{\{1, 2, 3, 4\}}^\text{all-rounder} \}}_\text{team 2} 
\end{equation*}
It will be crucial to see that a suitable number of teams $n$ can always be determined.

\begin{theorem}
    EDUCATIONAL TEAM FORMATION is \textbf{NP}-complete.
\end{theorem}

\begin{proof}
    It remains to give the details of the reduction from SET COVER sketched above.
    Let $(U, \mathcal{S}, k)$ be a SET COVER instance and set $m_0 = \#\mathcal{S}$ for an initial number of students with personal skill sets $S_a$ as given in $\mathcal{S}$.
    In the corresponding EDU-TF instance, we additionally fix $S = U$ as the set of all skills, $c = \#U$ as the minimum skill coverage, and $k_{max} = k$ is the upper bound for all team sizes, while we allow $k_{min} = 1$.
    Then we add another $(n-1)$-many students with $S_a = U$ and define them hereby as all-rounders (while $n$ is not fixed yet).
   
    The remaining step is to determine a suitable number of teams $n$ (and hence $(n-1)$ all-rounders) such that the total number of students  $(m_0+n-1)$ fit into $n$ teams while respecting the team-size constraint.
    This yields the following requirement for the choice of $n$:
    $$
    1 = k_{min} \leq \frac{m_0 + n - 1}{n} \leq k_{max} = k
    \;\;\; \iff \;\;\;
    1 \leq \frac{m_0 - 1}{k - 1} \le n
    $$
    Recall that $m_0 = \#\mathcal{S} \geq k \geq 2$, so we choose $n = \ceil{\frac{m_0-1}{k - 1}}$ which completes the construction of our EDU-TF instance.
    Observe that all steps of the construction are computable in polynomial time for any given SET COVER instance. 
    
    Now suppose that $(U, \mathcal{S}, k)$ is a yes-instance for SET COVER with solution $\mathcal{C}$.
    A solution to the corresponding EDU-TF instance looks as follows:
    One team consists of the $\#\mathcal{C}\leq k$ many students with skill sets in $\mathcal{C}$, while each of the remaining $n-1$ teams has exactly one all-rounder, hence the team-skill constraint is satisfied for all teams.
    The remaining $(\#\mathcal{S}-\#\mathcal{C})$ students can be arbitrarily assigned to the $n$ teams within the team-size constraint since their remaining capacity is
    \begin{align*}
    & \underbrace{(k-\#\mathcal{C})}_\text{team 1} + \underbrace{(n-1)(k-1)}_\text{teams $2,\ldots,n$} \\
    ={} & n(k-1)+1-\#\mathcal{C} \\
    ={} & \ceil{\frac{\#\mathcal{S}-1}{k - 1}}(k-1)+1-\#\mathcal{C} \\
    \geq{} & \#\mathcal{S}-\#\mathcal{C}.
    \end{align*}

    Conversely, note that in every solution to a yes-instance of EDU-TF constructed in the reduction
    there is at least one team with none of the additionally introduced all-rounders. So there exists a $(\leq k)$-choice from $\mathcal{S}$ that yields the universe $U$.
\end{proof}

The proof also shows that even the special case of EDU-TF is \textbf{NP}-complete where the lower team-size bound is trivial ($k_{min}=1$) and each team needs to cover all skills ($c=\# S$).

\subsection{Optimizing student preferences}

As pointed out earlier, we additionally take student preferences into account and extend the feasibility problem to an optimization problem.
For notational convenience, let
\begin{equation*}
    M = \{ (a, b) \in [m]^2 ~|~ a \not = b \land \exists T \in \mathcal{T} : a, b \in T \}
\end{equation*}
denote all pairs of students that are on the same team in some given feasible solution $\mathcal{T}$.
Furthermore, a given preference $p_{a, b}\in [-d, d]$ is called \emph{realized} in solution $\mathcal{T}$ if students $a$ and $b$ are assigned to the same team, i.e., $(a, b) \in M$.
Given a feasible solution $\mathcal{T}$ for EDU-TF we identify the following three objectives that cover many scenarios in practice:
\begin{itemize}
    \item \textbf{O1 -- Maximize teammate preference satisfaction.}
    Maximize the sum of all realized preference values, i.e., $\sum_{(a, b) \in M} p_{a, b}$.
    \item \textbf{O2 -- Minimize teammate preference dissatisfaction.} 
    Maximize the smallest realized preference, i.e., $\min_{(a, b) \in M} p_{a, b}$.
    We note here that in the edge case that $M = \emptyset$, i. e. if every student is on their own team, the objective value of this objective is defined as $\max p_{a, b} + 1$ so that singleton-teams are considered optimal.
    \item \textbf{O3$_{p_0}^{+/-}$ -- Maximize/minimize specific preferences.}
    Maximize ($+$) or minimize ($-$) the number of realized preferences with a value of $p_0$, i.e., $\#\{(a, b) \in M ~|~ p_{a, b} = p_0 \}$.
\end{itemize}
These objectives aim to result in teams that correspond to the two main aspects that can be seen in the student preferences: Students wanting to be teamed up with some students (O1, O3$_{p_0}^{+}$) and not wanting to be teamed up with other students (O1, O2, O3$_{p_0}^{-}$).
For $O\in \{\text{O1},\text{O2},\text{O3}\}$ we define optimization problems EDU-TF($O$) as follows.

\begin{problem}
    EDUCATIONAL TEAM FORMATION($O$) (EDU-TF($O$))\\
    Given the team, skill and preference data, provide a feasible EDUCATIONAL TEAM FORMATION solution that is optimal w.r.t. $O$, or report that no feasible solution exists.
\end{problem}

Since all of the above objective functions can be computed in polynomial time, each of the optimization problems EDU-TF($O$) is an \textbf{NP}-optimization problem.
Furthermore, they all extend the feasibility problem EDU-TF, so each EDU-TF($O$) is \textbf{NP}-hard (for basics notions of \textbf{NP}-optimization problems we refer to \citet{au03}).
In the following, we will consider multiple objectives in a hierarchical way and denote the respective optimization problems as EDU-TF($O_1, O_2, \ldots$).
Here each objective $O_i$ restricts the search space for $O_{i+1}$ to all solutions that are optimal w.r.t. $O_i$.
For different notions of multi-objective optimization and their mutual relations regarding hardness and separation we refer to \cite{10.1007/978-3-642-13962-8_20}.
Due to this hierarchical approach, it is possible to create different strategies in accordance with the requirements of the teacher or the students by choosing and ordering the objectives accordingly.
It also allows the value of the objective function in the different stages of the optimization process to still be easily interpreted.

\section{Solving the EDUCATIONAL TEAM FORMATION problem with a modular integer linear program}
\label{sec:ilp}

We solve the EDUCATIONAL TEAM FORMATION problem by means of an \textit{integer linear program} (ILP).
This ILP is constructed in two steps.
At first, only the feasibility problem EDU-TF is considered.
Its ILP will serve as the base model which can be extended by different objective functions.
In the second step, the objectives O1, O2, and O3 are modeled individually.
The result is a modular ILP that can be used to solve the optimization problem EDU-TF($O$) for any of the objectives O1, O2, and O3.
In addition, ILP solvers typically support multiple objective functions by solving them hierarchically, i.e., once the optimization of one objective finishes (e.g., by finding an optimal solution or by exhausting a configurable time limit) with a final objective value $v$, the solutions of all subsequent optimizations are constrained to have at least an objective value of $v$ w.r.t. that objective.
As such, the extended optimization problem EDU-TF($O_1, O_2, \ldots$) can also be solved for any number of objectives.
This allows us to explore multiple strategies by combining and optimizing different objectives in a hierarchical manner.

The following ILPs incorporate the problem's input parameters as presented in section \ref{sec:problem-data}.
To ease readability, we make use of indicator, minimum and logical constraints, as these are typically directly supported by ILP solvers.
Additionally, any constraints on the domain of a variable won't be explicitly listed among the other constraints.

\subsection{The base model}

The base model only deals with finding a feasible solution for the EDU-TF problem.
It uses the following variables:
\begin{alignat*}{2}
    x_{a, j} & \in \{0, 1\} \quad && \text{indicates whether student $a$ is assigned to team $j$.} \\
    y_{j, i} & \in \mathbb{N} \quad && \text{counts how many members of team $j$ cover skill $i$.} \\
    z_{j, i} & \in \{0, 1\} \quad && \text{indicates whether team $j$ covers skill $i$.}
\end{alignat*}
Using these variables and the EDU-TF problem's input parameters, the base ILP model includes the following constraints:
\begin{alignat}{2}
     \sum_{j \in [n]} x_{a, j} & = 1 & \quad \forall a & \in [m] \\
     \sum_{a \in [m]} x_{a, j} & \geq k_{min} & \quad \forall j & \in [n] \\
     \sum_{a \in [m]} x_{a, j} & \leq k_{max} & \quad \forall j & \in [n] \\
     y_{j, i} & = \sum_{a \in [m]: i \in S_a} x_{a, j} & \quad \forall j & \in [n], i \in [l] \\
     z_{j, i} & = \min(1, y_{j, i}) & \quad \forall j & \in [n], i \in [l] \\
     \sum_{i \in [l]} z_{j, i} & \geq c & \quad \forall j & \in [n]
\end{alignat}
Constraint (3) ensures that each student $a$ is assigned to exactly one team while (4) and (5) keep each team's size within the specified bounds.
The team skill constraint is modeled by constraints (6) through (8).
Constraints (6) counts the members of team $j$ covering skill $i$.
This count is transformed into the binary indicator $z_{j, i}$ in constraint (7), which is then used in constraint (8) to ensure each team covers a minimum of $c$ skills.

\subsection{Extensions to the base model for each objective}

We now expand the base model to take our objecitves into account.

\paragraph{O1 -- Maximize teammate preference satisfaction}
Objective O1 maximizes the sum of realized preferences.
In order to track which preferences are realized, the following variables are used:
\begin{alignat*}{2}
    t_{a, b, j} & \in \{0, 1\} \quad && \text{indicates whether students $a$ and $b$} \\
    &&& \text{are both assigned to team $j$.}
\end{alignat*}
The objective is then modeled as follows:
\begin{alignat}{3}
     \text{max} \quad && \sum_{a,b \in [m], j \in [l]} t_{a, b, j} p_{a, b} \span \\
     \text{s.t.} \quad && t_{a, b, j} & = x_{a, j} \land x_{b, j} & \quad \forall a, b & \in [m], j \in [n]
\end{alignat}
Constraint (10) ensures that $t_{a, b, j}$ is set to 1 if and only if both students $a$ and $b$ are assigned to team $j$.
Objective (9) models all realized preferences.
Note that the base model's constraints ensure that for each student pair $(a, b)$ exactly one of the indicator variables $t_{a, b, j}$ is set to 1.

\paragraph{O2 -- Minimize teammate preference dissatisfaction}
With objective O2, we seek to avoid assigning two students who do not wish to work together to the same team.
This requires the following variables:
\begin{alignat*}{2}
    q_{a, b} & \in \{0, 1\} \quad && \text{indicates whether students $a$ and $b$} \\
    &&& \text{are assigned to the same team.} \\
    r & \in \mathbb{R} \quad && \text{models the smallest realized preference in the solution.}
\end{alignat*}
The base model is then extended as follows:
\begin{alignat}{3}
    \text{max} \quad && r \span \\
    \text{s.t.} \quad && q_{a, b} = 1 & \implies x_{a, j} = x_{b, j} & \quad \forall a, b & \in [m], j \in [n] \\
    && q_{a, b} = 0 & \implies x_{a, j} + x_{b, j} \leq 1 & \quad \forall a, b & \in [m], j \in [n] \\
    && q_{a, b} = 1 & \implies r \leq p_{a, b} & \quad \forall a, b & \in [m]
\end{alignat}
Constraints (12) and (13) are indicator constraints and ensure that the indicator $q_{a, b}$ is set to 1 if and only if students $a$ and $b$ are assigned to the same team.
Constraint (14) bounds $r$ to the preference $p_{a, b}$ between $a$ and $b$ if they are assigned to the same team.
Hence, maximizing $r$ ensures that optimal solutions avoid low-preference students on the same team whenever feasible.

\paragraph{O3$_{p_0}^{+/-}$ -- Maximize/minimize specific preferences}
Objective O3 maximizes ($+$) or minimizes ($-$) the number of realized  preferences with a value of $p_0$. 
Let $A = \{(a, b) \in [m]^2 ~|~ p_{a, b} = p_0\}$ be the set of student pairs with such preferences.
Then only for all $(a, b) \in A$ the following variable is introduced:
\begin{alignat*}{2}
    t_{a, b, j} & \in \{0, 1\} \quad && \text{indicates whether students $a$ and $b$} \\
    &&& \text{are both assigned to team $j$.}
\end{alignat*}
The base model extension is as follows:
\begin{alignat}{3}
    \text{max/min} \quad && \sum_{(a,b) \in A, j \in [l]} t_{a, b, j} \span \\
    \text{s.t.} \quad && t_{a, b, j} & = x_{a, j} \land x_{b, j} & \quad \forall (a, b) & \in A, j \in [n]
\end{alignat}
Objective (15) maximizes or minimizes the number of student pairs of interest assigned to the same team.
Constraint (16) is the same as (10).

\subsection{Implementation details}

We used the Gurobi solver version 10.0.1 for the implementation.
This heavily influenced how we modeled the indicator variable for two students being assigned to the same team.
Preliminary tests showed that using AND constraints for objectives O1 and O3 lead to shorter optimization times compared to using indicator constraints, while the opposite was true for objective O2.
Apart from this, we did not further experiment with different modeling strategies.

When extending the base model with multiple objectives, the variables and constraints of each objective are always added in full without considering potential redundancies.
E.g., the indicator variables $t_{a, b, j}$ of objective O1 are not reused for objective O3.
Instead, if O3 is used in conjunction with O1, the variables $t_{a, b, j}$ are added again with constraint (16).
While this results in an overall larger model, it also facilitates a more straight-forward implementation, allowing for an easier adoption of our approach.
Additionally, this redundancy is resolved during a presolve phase performed by Gurobi.

\section{Experimental results}
\label{sec:results}

The core considerations of the proposed approach are its applicability to a variety of real-world scenarios and its adaptability to different teacher and student needs.
For its evaluation, we applied our model to multiple project-based courses of varying sizes.
These courses required different team sizes and employed different methods for collecting student preferences, ranging from explicit teammate requests to using the similarities between the students' profiles.
A multitude of different objective combinations, which we refer to as \emph{strategies}, were used for the optimization of the team composition for each course, highlighting the flexibility of our approach.

The evaluated strategies are analyzed and compared against the manual, teacher-assigned strategy, using optimization runtime and the number of realized student preferences as evaluation metrics.
Other metrics, such as final team performance or grade outcomes, were not considered, primarily due to the ethical considerations present in educational settings.
Specifically, conducting a controlled study with different team formation methods for separate student groups (e.g., manual teacher assignment for one part of the class versus algorithmic assignment for another) was deemed unfair for students, as it could result in different academic outcomes and unequal learning experiences.
Similarly, a longitudinal study comparing historical and current team formation methods was infeasible across all courses.
Our emphasis was therefore on examining the algorithm’s performance across varied datasets and strategies, rather than assessing its impact on academic performance in only one or two courses.

The rest of this section is structured as follows.
Subsection \ref{sec:datasets} presents the different datasets that were derived from the individual courses.
The evaluated strategies are described in subsection \ref{sec:strategies}. Finally, subsection \ref{sec:observations} concludes with detailed observations and a comprehensive comparison of the different strategies for each dataset.
A set of key takeaways, along with a broader discussion of the contributions and limitations of our approach, is included in the next section, \ref{sec:discussion}.

\subsection{Datasets}
\label{sec:datasets}

We worked with datasets from eight university-level, project-based courses on computing education, categorized into nine distinct datasets (see Table~\ref{tab:overview-datasets}).
The courses ranged from small-sized courses consisting of $m=15$ students to mid-sized courses with $m=79$ students.
They also varied in terms of their teacher requirements regarding the team size, ranging from student pairs (with exceptionally-allowed single-student teams in an odd-sized cohort) to teams of up to $k_{max}=6$ students, resulting in the formation of $n=7$ to $n=27$ teams.
Most courses differentiated between $\#S=2$ and $\#S=6$ skills with one course defining $\#S=14$ distinct skills.
As long as a feasible solution exists, all skills had to be covered ($c=\#S$), otherwise a smaller minimum skill coverage was chosen.
Skills were gathered by first asking the students to assess their perceived skill level on a numerical scale before applying a threshold to indicate whether a student is considered to cover a skill.
The thresholds were defined per skill by the instructor of the course such that each skill was covered by a sufficient number of students.

Finally, the courses employed different methods for gathering the student's preferences.
In most cases, every student was able to explicitly provide preferences for a small subset of their classmates, resulting in a sparse preference matrix (i.e., most entries are 0).
In order to obtain a dense preference matrix, many courses additionally required the students to create a personal profile based on a set of predefined questions.
The similarity of these profiles could then be used as a proxy to derive preferences in the cases where no explicit preference was provided, enabling the creation of dense preference matrices.
These profiles are the differentiating factor between the datasets D8.1 and D8.2, which originate from the same course.
D8.1 includes the profile based preferences while D8.2 does not, permitting us to experiment with the algorithm’s performance with a dense and sparse preference matrix respectively in the presence of a large student cohort.

The process of gathering the students' profile attributes and their explicit preferences are described in the following.

\subsubsection{Explicit preference extraction}
The explicit preferences can be categorized into ``weak'' and ``strong'' preferences.
Weak preferences were gathered through a process referred to as ``team dating'' \citep{umbelino2021prototeams,lykourentzou2017team}; a method allowing small-group interactions (in this case pairs) to enable participants to obtain, in a short time frame, direct experience over multiple elements of their potential collaboration and assess whether they could collaborate efficiently.
The team dating sessions were organized for the courses that this was possible due to course restrictions (i.e., availability of teaching slots and/or class size), and usually lasted for one hour during which each student was able to discuss with four potential teammates.
At the end of the process, students gave their preference rating on a Likert scale from 1 to 5 to the fellow students they talked to.
These preferences were then remapped to the interval $[-2,2]$.

In addition to these weak preferences, some instructors allowed their students to directly provide a list of students they did want and did not want to work with.
These preferences are considered ``strong'' as they stem from previous experiences and should be prioritized over the preferences resulting from the students' first interaction during the team dating session.
As such, these strong preferences were assigned the values 4 and $-4$.

\subsubsection{Profile attribute extraction}
The profiles were gathered through a questionnaire, inquiring about the students' a) time management preference (i.e., whether they prefer to start working on the group project as early as possible or closer to the deadline), b) intended effort in the course and c) preference for synchronous or asynchronous collaboration.
While the latter question only allowed for a binary choice, the time management preference was measured on a Likert scale from 1 to 5 and the intended effort utilized a scale from 1 to 3.
For each course, the following steps were then performed to calculate preferences from these profiles:

\begin{enumerate}
    \item The profiles were normalized by rescaling each attribute to the interval $[0, 1]$.
    This was done to ensure that all attributes are weighed equally in the similarity calculation.
    The normalization was done based solely on the ratings actually present in the original profiles.
    \item We used the euclidean distance between the profile vectors to calculate their similarity.
    In order to obtain naturally interpretable similarities, the calculated distances were rescaled to $[0, 1]$ (with 1 being the maximum distance among all profiles).
    \item The range of profile similarities was discretized into a fixed number of evenly sized buckets so that similarities could be mapped to preferences according to their bucket.
\end{enumerate}

In case a course did not allow the students to supply weak preferences through a team dating session, the profiles were used as a replacement for the weak preferences.
As such, the profile similarities were mapped to the preference interval $[-2, 2]$.
When there were already explicit weak preferences present, the profiles were used to supplement them.
Because we wanted to prioritize these explicit preferences over the similarity based ones, the profile similarities were mapped to the smaller preference interval $[-1, 1]$.

\begin{table}
    \centering
    \begin{tabular}{lccccccccccccc}
        \toprule
                      & \multicolumn{5}{c}{EDU-TF input parameters} & \multicolumn{7}{c}{Preferences} \\ \cmidrule(lr){2-6} \cmidrule(lr){7-13}
                      & $m$ & $n$ & $[k_{min}, k_{max}]$ & $\#S$ & $c$ & 4   & 2    & 1    & 0    & $-1$ & $-2$ & $-4$ & Profiles \\
        \midrule
        \textbf{D1}   & 15  & 8   & $[1, 2]$             & 2     & 2   & --- & 8    & 70   & 76   & 56   & 0    & ---  & yes      \\ 
        \textbf{D2}   & 26  & 13  & $[2, 2]$             & 3     & 3   & --- & 110  & 152  & 88   & 266  & 34   & ---  & yes      \\ 
        \textbf{D3}   & 32  & 16  & $[2, 2]$             & 3     & 3   & --- & 26   & 49   & 899  & 13   & 5    & ---  & no       \\ 
        \textbf{D4}   & 32  & 7   & $[4, 5]$             & 3     & 3   & 37  & 207  & 253  & 55   & 319  & 121  & ---  & yes      \\ 
        \textbf{D5}   & 37  & 10  & $[3, 4]$             & 14    & 3   & 6   & 110  & 285  & 525  & 344  & 52   & 10   & yes      \\ 
        \textbf{D6}   & 60  & 12  & $[4, 6]$             & 6     & 4   & 138 & 96   & 81   & 3018 & 94   & 113  & ---  & no       \\ 
        \textbf{D7}   & 69  & 18  & $[3, 4]$             & 4     & 4   & 4   & 43   & 1490 & 2272 & 881  & 2    & ---  & yes      \\ 
        \textbf{D8.1} & 79  & 27  & $[2, 3]$             & 3     & 3   & 23  & 25   & 1923 & 2961 & 1223 & 7    & ---  & yes      \\ 
        \textbf{D8.2} & 79  & 27  & $[2, 3]$             & 3     & 3   & 23  & 25   & 50   & 6036 & 21   & 7    & ---  & no       \\ 
        \bottomrule
    \end{tabular}
    \caption{Overview of the used courses and datasets. The input parameters are as follows: $m$ -- number of students; $n$ -- number of teams; $k_{min}, k_{max}$ -- team size bounds; $\#S$ -- number of skills; $c$ -- number of skills to be covered.}
    \label{tab:overview-datasets}
\end{table}

\subsection{Strategies}
\label{sec:strategies}

The modular structure of our model facilitates the hierarchical optimization of multiple objectives as explained above.
Going forward, we call an objective or a combination of objectives a \emph{strategy} for solving the EDUCATIONAL TEAM FORMATION problem from a practical point of view, as the objectives (and their order) express a constructive process for how to create a team composition that is considered ``good'' by a teacher or student.
By choosing the objectives and their order accordingly, a teacher can create or adapt a strategy specifically for the team formation requirements of a course.
For example, if minimizing the dissatisfaction in the teams is prioritized, a teacher could use objective O2 to minimize the dissatisfaction first, followed by O1 to maximize the satisfaction while keeping the dissatisfaction minimal.
We call this strategy EDU-TF(O2, O1) (i.e., the same as the corresponding formal problem), as it uses O2 and O1 with the objectives being optimized from left to right.

Table \ref{tab:overview-strategies} gives an overview over all strategies we explored.
For every dataset we calculated solutions with every strategy presented in the table.
The strategies were slightly adapted to datasets that lack some preference values.
In these cases, the use of objective O3 to maximize or minimize the number of realized preferences with that value was skipped, e.g., strategy S1.2 for dataset D1 is EDU-TF(O3$_{-2}^-$, O3$_{-1}^-$) instead of EDU-TF(O3$_{-4}^-$, O3$_{-2}^-$, O3$_{-1}^-$), because no preferences with value $-4$ were given.
This way, the general idea of each strategy remained unaltered.
For datasets without strong preferences (preferences with a value of $4$ or $-4$), this results in strategies S2.1 and S4.1 being the same.
The strategies we used can be divided into the following four groups:

\begin{table}
    \centering
    \begin{tabular}{llp{180pt}}
        \toprule
        ID & Strategy & Description \\
        \midrule
        \textbf{S1.1} & EDU-TF(O2) & 
        \multirow{2}{=}{
        Minimize dissatisfaction} \\
        \textbf{S1.2} & EDU-TF(O3$_{-4}^-$, O3$_{-2}^-$, O3$_{-1}^-$) & \\
        \midrule
        \textbf{S2.1} & EDU-TF(O1) & \multirow{2}{=}{Maximize satisfaction} \\
        \textbf{S2.2} & EDU-TF(O3$_{4}^+$, O3$_{2}^+$, O3$_{1}^+$) & \\
        \midrule
        \textbf{S3.1} & EDU-TF(O2, O1) &  \multirow{4}{=}{Minimize dissatisfaction, then maximize satisfaction} \\
        \textbf{S3.2} & EDU-TF(O2, O3$_{4}^+$, O3$_{2}^+$, O3$_{1}^+$) & \\
        \textbf{S3.3} & EDU-TF(O3$_{-4}^-$, O3$_{-2}^-$, O3$_{-1}^-$, O1) & \\
        \textbf{S3.4} & EDU-TF(O3$_{-4}^-$, O3$_{-2}^-$, O3$_{-1}^-$, O3$_{4}^+$, O3$_{2}^+$, O3$_{1}^+$) & \\
        \midrule
        \textbf{S4.1} & EDU-TF(O3$_{-4}^-$, O3$_{4}^+$, O1) &  \multirow{2}{=}{Optimize strong preferences first}\\ 
        \textbf{S4.2} & EDU-TF(O3$_{-4}^-$, O3$_{4}^+$, O3$_{-2}^-$, O3$_{2}^+$, O3$_{-1}^-$, O3$_{1}^+$) & \\
        \bottomrule
    \end{tabular}
    \caption{Overview over the used strategies}
    \label{tab:overview-strategies}
\end{table}

\paragraph{Minimize dissatisfaction}
The first group contains strategies that only try to avoid teaming up students with a low-valued preference between them without considering the high-valued preferences.
This can either be done by just using objective O2 (S1.1), or objective O3 for each negative preference to explicitly minimize the realization of those preferences in order from $-4$ to $-1$ (S1.2).
S1.2 can easily be adapted to other potential datasets, in which a different scale of preferences is used, by adding, changing, or leaving out variations of O3.
S1.1 can be used for every dataset without adaptation.
There are two key differences between the two strategies.
While S1.2 only minimizes the preferences that are explicitly given, S1.1 can result in the elimination of neutral or small positive preferences should feasible solutions like this exist.
On the other hand, if some negative preferences are unavoidable, S1.2 will still seek to minimize the number of realizations of these unavoidable and undesirable preferences, whereas S1.1 will potentially produce solutions that realize a lot of these preferences since they cannot be completely eliminated.

\paragraph{Maximize satisfaction}
The strategies in the second group maximize the student satisfaction by teaming up students with high-valued preferences between them.
S2.2 does this by explicitly maximizing the positive-valued preferences, starting with the highest.
Just like S1.2, strategy S2.2 can easily be adapted to other preference scales.
S2.2 only considers the positive-valued preferences given in the strategy.
In contrast, using objective O1, strategy S2.1 maximizes the satisfaction by optimizing the overall sum of all realized preferences.
Because of this, S2.1 is the simplest strategy that takes all preferences into account.
However, it is possible that S2.1 realizes less positive preferences than S2.2 if realizing more positive preferences results in also realizing more negative ones, resulting in a smaller overall sum.

\paragraph{Minimize dissatisfaction, then maximize satisfaction}
In group three, the strategies from group one and two are combined to create strategies that consider (in most cases) all preferences, while prioritizing the minimization of the student dissatisfaction over the student satisfaction.
This results in the four strategies seen in Table \ref{tab:overview-strategies}. 
These strategies almost always consider all preferences in their calculations.
The exception is S3.2 if not all negative preferences can be eliminated.
In these cases, the negative preferences that are not completely eliminated are not minimized.

\paragraph{Prioritize strong preferences}
Similarly to group three, the strategies in group four take all preferences into account, with the strong preferences being prioritized.
Therefore, both strategies start by explicitly minimizing preferences with a value of $-4$ and maximizing preferences with a value of $4$.
After that, S4.1 optimizes the other preferences with objective O1.
In comparison, S4.2 explicitly minimizes or maximizes the other preferences using O3 as well, in descending order of their absolute value, always starting with the negative preference.

Besides the strategies we used, it is also possible to create more strategies, e.g., if the maximization of student satisfaction should be prioritized, the combination of strategies in group three can be done the other way around.

\subsection{Strategy comparison}
\label{sec:observations}

We evaluate each of the ten proposed strategies on every dataset using the following metrics: 1) the quality of the resulting solution as defined by the sum of all realized preferences, and 2) the runtime of the optimization.
In our analysis, we also highlight more nuanced observations based on the individual preference values.

Each strategy was applied twice using two different methods to limit the runtime of the computation.
The first method used a time limit of one hour for the entire optimization, i.e., each individual objective was solved optimally before progressing with the next objective with the entire optimization being interrupted after one hour.
This one-hour time limit allows for many strategies to produce an optimal solution while still limiting the overall time required for running these tests.
The second method used a time limit of 15 minutes divided equally among all objectives.
Thus, if a strategy contains three objectives, each objective is assigned a time limit of 5 minutes.
This approach is motivated by the need to quickly create teams in volatile settings.
By limiting the runtime for each objective, we expected the resulting solution to still conform to the priorities expressed by the overall strategy and to be of sufficient quality despite the small time frame.
Note that if the optimization of an objective finished earlier (i.e., if an optimal solution for that objective was found), the remaining time was \emph{not} redistributed among the remaining objectives.

Each dataset includes a hand-crafted solution, which serves as a baseline for assessing the effectiveness of the strategies.
It is of note, however, that not all manual solutions were created using the identical problem instances as the ones the strategies were applied to.
The preference matrix of most datasets combines the explicitly stated preferences by the students with preferences derived from their profiles.
In these instances, however, the instructors responsible for creating the manual team assignments did not have access to the entire preference matrix and instead primarily focused on the explicit preferences.
Nevertheless, they still had access to the raw profiles as well and potentially used them to aid in their decision making.
As such, we still deem it appropriate to compare the algorithmic solutions to the manual ones in these cases.

The results for the one-hour and fifteen-minute time limit are presented in table \ref{tab:exp_results}.
The manual solutions that were created without the profile based preferences are marked with an asterisk (*).
An asterisk is additionally used to highlight strategy S4.1 for the datasets D1, D2 and D3.
This is due to the preference values 4 and $-4$ not being present.
As a result, strategy S4.1 is the same as S2.1, since S4.1 only differs in explicitly optimizing for these preferences first.

\paragraph{General observations}
Strategy S2.1, which explicitly maximizes the sum of all realized preferences without any other objective, naturally resulted in the solution with the highest quality for each dataset compared to the other solutions with the same runtime restrictions.
In some datasets there is a trade-off between achieving the highest quality and optimizing some of the other objectives. 
For instance, strategies like S3.3, which do maximize the sum of the realized preferences in the end, could not produce solutions with the same overall quality, even if they were solved optimally.
Nevertheless, in most cases the strategies S2.2 - S4.2 resulted in solutions of similar quality compared to S2.1.
Almost all of the solutions for these strategies and both time limits outperform the manual solutions given for each dataset.
The only ones that usually do not outperform the quality of the manual solutions are the solutions for strategies S1.1 and S1.2.
They primarily minimize the realization of negative preferences and, therefore, do not result in a high quality solutions in most cases.
They can be used, however, as a baseline for which negative preferences can be avoided, although they should not be used to form teams.
Finally, comparing the runtime of these two strategies highlights a significant difference.
S1.1 required a noticeably longer runtime that S1.2 (except for the smallest datasets), with the difference becoming more extreme for larger datasets.


\paragraph{D1, D2, D3}
For the three smallest datasets, 15 to 32 students were partitioned into teams of size two.
The runtime in most cases was just a few seconds.
In cases in which the whole calculations took considerably longer (e.g., S3.2 and S3.4 for D2), the optimal solution was found in under 30 seconds.
The rest of the time was used to prove that the solution is optimal.
Because all calculations finished within the time limit, all solutions are optimal regarding the different combinations of objectives.
For these datasets, strategies S2.1 - S4.2 all resulted in solutions with the same overall quality and only slight differences in the numbers of realized preferences, while the runtimes display more variance, especially in the case of D2.
Compared to most of the other calculations, strategies S3.1 and S3.2 needed more time to find the best solution and to prove their optimality.
Due to the overall short runtime, there is no difference in runtimes or solutions between the one-hour and fifteen-minute time limit calculations for D1 and D3.
For D2, however, the runtime for the fifteen-minute calculations was shorter in some cases, while still resulting in the same solutions.

\paragraph{D4, D6}
With the exception of S3.1 for D4, all calculations for the datasets D4 and D6 finished within the one-hour time limit and, therefore, were solved optimally.
The solutions for D4 were found in under one minute, and there is no difference between the one-hour and the fifteen-minute calculations.
For D6 the runtimes were more diverse, e.g., S3.3 finished in under one minute, while the best solution for S3.1 was found in about 12 minutes and was shown to be optimal after 34 minutes.
Here, equally dividing the fifteen-minute time limit over all objectives could improve the runtime, especially if it took a long time to prove that a solution is optimal, during which the fifteen-minute calculations often reached the time limit.
Despite not solving all objectives provably optimal with a time limit of 15 minutes, the resulting solutions are still optimal regarding the individual objectives, as they are identical to the optimal solutions found during the longer calculations.
For both datasets strategies S3.1 and S3.2, which start with objective O2, took by far the longest to find the best solution, and they had the longest overall runtime.
The most noteworthy case is S3.1 for D4, which did not finish within one hour.

Taking a closer look at the solutions for both datasets, a trade-off between the different objectives can be seen, which is more extreme in D4.
The highest quality, which is reached with strategy S2.1, cannot be reached anymore if, before maximizing the sum of realized preferences, the negative preferences are minimized (S3.3) or the number of realized preferences with a value of 4 is maximized (S4.1).

\paragraph{D5}
D5 is the only dataset that does include preferences with a value of $-4$ and that defines the largest set of skills, of which only three need to be covered per team.
This dataset exhibits the most extreme difference between the overall runtime and the time at which the final improvement to a solution is found.
These improvements occurred after at most one and a half minutes, but almost all calculations timed out after one hour with the exception of S3.1 and S4.1, which finished after about 20 seconds.
Applying strategies S3.1 and S3.2 required more time to find the best solution or to complete than their corresponding strategies S3.3 and S3.4 that do not use objective O2, both in the one-hour and the fifteen-minute calculations.
Other than S3.2, the calculations that timed out solved all but the last objective of their strategy optimally. 
Nevertheless, their solutions are very similar to the optimal solutions for S2.1 and S4.1.
In the calculation for S3.2, the time limit was reached before the final objective was considered, which also resulted in this solution being slightly worse than the others.
The fifteen-minute calculations for all strategies resulted in the exact same solutions as the one-hour calculations except for S3.2.
In fact, the solution for S3.2 was improved compared to the longer time limit, because all objectives are considered in the fifteen-minute calculations due to limiting the runtime per objective.
Comparing the optimal solutions for S3.1 and S4.1, a small trade-off between the overall quality and the number of realized preferences with a value of 4 can be observed.
If all six of these preferences are realized, the highest possible quality of 149 can no longer be achieved.

\paragraph{D7, D8.1}
D7 and D8.1 are the two biggest datasets with 69 and 79 students respectively.
Apart from S1.1 and S1.2, all optimizations reached the time limit of one hour, with the final improvement happening fairly late into the optimization process, suggesting that further improvements are likely possible.
The timeout most frequently occurred during the optimization of the final objective.
The exception to this is S3.2.
With a time limit of 15 minutes, the optimizations of more objectives timed out.
The runtime of S3.1 and S3.2 for D8.1 significantly exceeded 15 minutes, which is due to a general presolve step performed by Gurobi that is not affected by the time limit for the individual objectives.

Most strategies produced solutions of similar quality that significantly outperform the manual solution.
Excluding S1.1 and S1.2, the strategies that resulted in noticeably worse solutions are S2.2 and S3.2 for D7 and S3.1 and S3.2 for D8.
The shorter time limit of 15 minutes overall resulted in worse solutions.
For D7, the most significant drop in quality can be observed, while the solutions for D8.1 are mostly comparable to the one-hour time limit.
The exception to this are S3.1 and S3.2, for which no significant improvements to the solution were found, resulting in no objective being solved optimally and the final solutions being worse than the manual baseline.
This suggests that strategies including O2 result in ILPs that are more difficult to solve in practice.
Finally, comparing S3.4 and S4.2, a small trade-off can be observed regarding the realization of more positive preferences and avoiding negative preferences.

\begin{figure}[ht]
    \centering
    \includegraphics[width=1\linewidth]{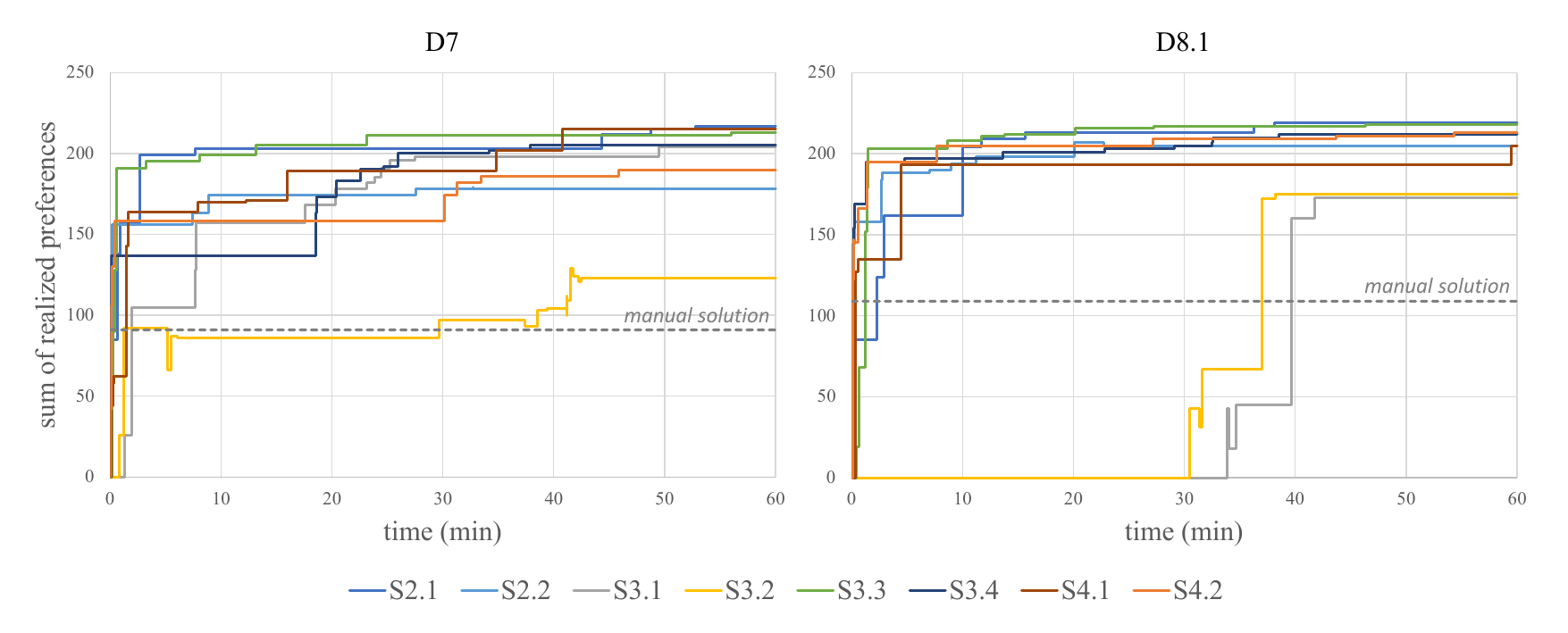}
    \caption{Evolution of the solution quality expressed as the sum of all realized preferences during the calculation of strategies S2.1 - S4.2 on datasets D7 and D8.1. The sums of realized preferences in the manual solutions are shown as horizontal dashed lines. We observe that the most substantial improvements of the sum occur in the beginning of the calculations and that the manual solution is outperformed by almost all strategies after about two minutes. For dataset D7 S3.2 is the worst performing strategy both in solution quality and in time needed to reach this quality. Strategy S3.2 and additionally S3.1 also perform the worst for dataset D8.1 with the biggest steps in solution quality happening after about 30 minutes instead of the beginning and the sum of realized preferences reached being smaller compared to the other strategies.}  
    \label{fig:calc_D7_D8.1}
\end{figure}

Since D7 and D8.1 are the only datasets for which the solutions improved up until the end of the one hour time limit, we further analyzed the evolution of the solution quality throughout the entire optimization in order to gain insight into when the most significant improvements occur.
Although most objectives don't directly optimize the sum of the realized preferences, we examined only this metric as it allows for a unified presentation across all objectives in a single graph and an improvement to any objective usually improves this sum as well.
The resulting graphs are shown in figure \ref{fig:calc_D7_D8.1}.
They highlight that the most substantial improvements happened at the start of the optimization, resulting in most strategies outperforming the manual solution after two minutes.
The strategies S3.1 and S3.2, however, do not align with this general observation.
When they are applied to D8.1, the first solution is found after around half an hour during the optimization of their common first objective O2.
This further underlines that ILPs incorporating O2 are likely more difficult to solve, making them less practical.

\paragraph{D8.2}
The datasets D8.1 and D8.2 consist of the same set of students and their explicitly given preferences.
While the preference matrix for D8.1 additionally contains preferences calculated from student profiles, the preference matrix for D8.2 only consists of the explicit preferences.
As a result, the matrix is much more sparse.
Compared to D8.1, for which most calculations timed out after one hour, nearly all strategies applied to D8.2 resulted in optimal solutions in at most 80 seconds.
The overall faster runtimes in comparison to D8.1 indicate that for larger datasets it could be beneficial to use sparse preference matrices in order to get optimal solutions using only the explicitly given preferences in a shorter amount of time.
Once again, strategies S3.1 and S3.2 performed the worst both in runtime and in solution quality.
For both of them, only the first objective O2 could be solved optimally in the time limit and the resulting solutions are one of the few exceptions that neither outperform the manual solution in the number of positive preferences being realized nor regarding the overall solution quality.
Other than S2.2, the rest of the strategies result in the same solution without negative preferences and an overall quality of 136.
The solution for S2.2 realizes more positive preferences at the cost of also realizing some negative preferences and achieving a slightly lower quality.
Comparing this solution to the solution for S3.4 shows that there is a trade-off between realizing positive preferences and avoiding negative preferences for this dataset.

\begingroup
\pagebreak

\setlength{\tabcolsep}{1.4pt}  
\newcommand{\xhspace}{\hspace{2\tabcolsep}}  
\renewcommand{\arraystretch}{0.96}  

\begin{longtable}[c]{lll@{\xhspace}|@{\xhspace}rrrrrrrr@{\xhspace}|@{\xhspace}ccl}
\toprule
    &              &               & \multicolumn{7}{c}{Realized preferences} &     &      $T$     &     $T^*$    &            \\
    & \multicolumn{2}{l|@{\xhspace}}{Strategy} &  4 &  2 &  1 &  0 &  $-1$ &  $-2$ & $-4$ & Sum & \multicolumn{2}{c}{[mm:ss]} & Objectives \\
\midrule
\endhead

\bottomrule
\endfoot
\endlastfoot

\multirow{22}{*}{\rotatebox{90}{\textbf{D1}}}
    & \multicolumn{2}{l|@{\xhspace}}{\textit{total}} & -- & \textit{8} & \textit{70} & \textit{76} & \textit{56} & \textit{0} & -- & & & & \\
    & \multicolumn{2}{l|@{\xhspace}}{manual*} & -- & 3 & 5 & 4 & 2 & \textbf{0} & -- & 9 & & & \\
\cmidrule(r){2-14}
& \multirow{10}{*}{\rotatebox{90}{1h}}
      & S1.1  & -- & 1 & \textbf{13} & 0 & \textbf{0} & \textbf{0} & -- & 15 & 00:01 & 00:01 & \textbf{O2} \\
    & & S1.2  & -- & 0 & 8 & 6 & \textbf{0} & \textbf{0} & -- & 8 & 00:01 & 00:01 & \textbf{O3$_{-2}^-$}, \textbf{O3$_{-1}^-$} \\
    & & S2.1  & -- & 3 & 11 & 0 & \textbf{0} & \textbf{0} & -- & \textbf{17} & 00:01 & 00:01 & \textbf{O1} \\
    & & S2.2  & -- & \textbf{4} & 9 & 1 & \textbf{0} & \textbf{0} & -- & \textbf{17} & 00:01 & 00:01 & \textbf{O3$_{2}^+$}, \textbf{O3$_{1}^+$} \\
    & & S3.1  & -- & 3 & 11 & 0 & \textbf{0} & \textbf{0} & -- & \textbf{17} & 00:02 & 00:01 & \textbf{O2}, \textbf{O1} \\
    & & S3.2  & -- & 3 & 11 & 0 & \textbf{0} & \textbf{0} & -- & \textbf{17} & 00:02 & 00:01 & \textbf{O2}, \textbf{O3$_{2}^+$}, \textbf{O3$_{1}^+$} \\
    & & S3.3  & -- & 3 & 11 & 0 & \textbf{0} & \textbf{0} & -- & \textbf{17} & 00:01 & 00:01 & \textbf{O3$_{-2}^-$}, \textbf{O3$_{-1}^-$, O1} \\
    & & S3.4  & -- & \textbf{4} & 9 & 1 & \textbf{0} & \textbf{0} & -- & \textbf{17} & 00:01 & 00:01 & \textbf{O3$_{-2}^-$}, \textbf{O3$_{-1}^-$}, \textbf{O3$_{2}^+$}, \textbf{O3$_{1}^+$} \\
    & & S4.1* & -- & 3 & 11 & 0 & \textbf{0} & \textbf{0} & -- & \textbf{17} & 00:01 & 00:01 & \textbf{O1} \\
    & & S4.2  & -- & \textbf{4} & 9 & 1 & \textbf{0} & \textbf{0} & -- & \textbf{17} & 00:01 & 00:01 & \textbf{O3$_{-2}^-$}, \textbf{O3$_{2}^+$}, \textbf{O3$_{-1}^-$}, \textbf{O3$_{1}^+$} \\
\cmidrule(r){2-14}
& \multirow{10}{*}{\rotatebox{90}{15min}}
      & S1.1  & -- & 1 & \textbf{13} & 0 & \textbf{0} & \textbf{0} & -- & 15 & 00:01 & -- & \textbf{O2} \\
    & & S1.2  & -- & 0 & 8 & 6 & \textbf{0} & \textbf{0} & -- & 8 & 00:01 & -- & \textbf{O3$_{-2}^-$}, \textbf{O3$_{-1}^-$} \\
    & & S2.1  & -- & 3 & 11 & 0 & \textbf{0} & \textbf{0} & -- & \textbf{17} & 00:01 & -- & \textbf{O1} \\
    & & S2.2  & -- & \textbf{4} & 9 & 1 & \textbf{0} & \textbf{0} & -- & \textbf{17} & 00:01 & -- & \textbf{O3$_{2}^+$}, \textbf{O3$_{1}^+$} \\
    & & S3.1  & -- & 3 & 11 & 0 & \textbf{0} & \textbf{0} & -- & \textbf{17} & 00:02 & -- & \textbf{O2}, \textbf{O1} \\
    & & S3.2  & -- & 3 & 11 & 0 & \textbf{0} & \textbf{0} & -- & \textbf{17} & 00:02 & -- & \textbf{O2}, \textbf{O3$_{2}^+$}, \textbf{O3$_{1}^+$} \\
    & & S3.3  & -- & 3 & 11 & 0 & \textbf{0} & \textbf{0} & -- & \textbf{17} & 00:01 & -- & \textbf{O3$_{-2}^-$}, \textbf{O3$_{-1}^-$, O1} \\
    & & S3.4  & -- & \textbf{4} & 9 & 1 & \textbf{0} & \textbf{0} & -- & \textbf{17} & 00:01 & -- & \textbf{O3$_{-2}^-$}, \textbf{O3$_{-1}^-$}, \textbf{O3$_{2}^+$}, \textbf{O3$_{1}^+$} \\
    & & S4.1* & -- & 3 & 11 & 0 & \textbf{0} & \textbf{0} & -- & \textbf{17} & 00:01 & -- & \textbf{O1} \\
    & & S4.2  & -- & \textbf{4} & 9 & 1 & \textbf{0} & \textbf{0} & -- & \textbf{17} & 00:01 & -- & \textbf{O3$_{-2}^-$}, \textbf{O3$_{2}^+$}, \textbf{O3$_{-1}^-$}, \textbf{O3$_{1}^+$} \\
\midrule

\multirow{22}{*}{\rotatebox{90}{\textbf{D2}}}
    & \multicolumn{2}{l|@{\xhspace}}{\textit{total}} & -- & \textit{110} & \textit{152} & \textit{88} & \textit{266} & \textit{34} & -- & & & & \\
    & \multicolumn{2}{l|@{\xhspace}}{manual*} & -- & 20 & 4 & 2 & \textbf{0} & \textbf{0} & -- & 44 & & & \\
\cmidrule(r){2-14}
& \multirow{10}{*}{\rotatebox{90}{1h}}
      & S1.1  & -- & 10 & \textbf{16} & 0 & \textbf{0} & \textbf{0} & -- & 36 & 00:11 & 00:02 & \textbf{O2} \\
    & & S1.2  & -- & 12 & 10 & 4 & \textbf{0} & \textbf{0} & -- & 34 & 00:01 & 00:01 & \textbf{O3$_{-2}^-$}, \textbf{O3$_{-1}^-$} \\
    & & S2.1  & -- & \textbf{24} & 2 & 0 & \textbf{0} & \textbf{0} & -- & \textbf{50} & 00:02 & 00:01 & \textbf{O1} \\
    & & S2.2  & -- & \textbf{24} & 2 & 0 & \textbf{0} & \textbf{0} & -- & \textbf{50} & 00:05 & 00:01 & \textbf{O3$_{2}^+$}, \textbf{O3$_{1}^+$} \\
    & & S3.1  & -- & \textbf{24} & 2 & 0 & \textbf{0} & \textbf{0} & -- & \textbf{50} & 01:28 & 00:19 & \textbf{O2}, \textbf{O1} \\
    & & S3.2  & -- & \textbf{24} & 2 & 0 & \textbf{0} & \textbf{0} & -- & \textbf{50} & 14:29 & 00:21 & \textbf{O2}, \textbf{O3$_{2}^+$}, \textbf{O3$_{1}^+$} \\
    & & S3.3  & -- & \textbf{24} & 2 & 0 & \textbf{0} & \textbf{0} & -- & \textbf{50} & 01:57 & 00:03 & \textbf{O3$_{-2}^-$}, \textbf{O3$_{-1}^-$, O1} \\
    & & S3.4  & -- & \textbf{24} & 2 & 0 & \textbf{0} & \textbf{0} & -- & \textbf{50} & 40:07 & 00:01 & \textbf{O3$_{-2}^-$}, \textbf{O3$_{-1}^-$}, \textbf{O3$_{2}^+$}, \textbf{O3$_{1}^+$} \\
    & & S4.1* & -- & \textbf{24} & 2 & 0 & \textbf{0} & \textbf{0} & -- & \textbf{50} & 00:02 & 00:01 & \textbf{O1} \\
    & & S4.2  & -- & \textbf{24} & 2 & 0 & \textbf{0} & \textbf{0} & -- & \textbf{50} & 04:19 & 00:01 & \textbf{O3$_{-2}^-$}, \textbf{O3$_{2}^+$}, \textbf{O3$_{-1}^-$}, \textbf{O3$_{1}^+$} \\
\cmidrule(r){2-14}
& \multirow{10}{*}{\rotatebox{90}{15min}}
      & S1.1  & -- & 10 & \textbf{16} & 0 & \textbf{0} & \textbf{0} & -- & 36 & 00:11 & -- & \textbf{O2} \\
    & & S1.2  & -- & 12 & 10 & 4 & \textbf{0} & \textbf{0} & -- & 34 & 00:01 & -- & \textbf{O3$_{-2}^-$}, \textbf{O3$_{-1}^-$} \\
    & & S2.1  & -- & \textbf{24} & 2 & 0 & \textbf{0} & \textbf{0} & -- & \textbf{50} & 00:02 & -- & \textbf{O1} \\
    & & S2.2  & -- & \textbf{24} & 2 & 0 & \textbf{0} & \textbf{0} & -- & \textbf{50} & 00:05 & -- & \textbf{O3$_{2}^+$}, \textbf{O3$_{1}^+$} \\
    & & S3.1  & -- & \textbf{24} & 2 & 0 & \textbf{0} & \textbf{0} & -- & \textbf{50} & 01:28 & -- & \textbf{O2}, \textbf{O1} \\
    & & S3.2  & -- & \textbf{24} & 2 & 0 & \textbf{0} & \textbf{0} & -- & \textbf{50} & 10:22 & -- & \textbf{O2}, o3$_{2}^+$, o3$_{1}^+$ \\
    & & S3.3  & -- & \textbf{24} & 2 & 0 & \textbf{0} & \textbf{0} & -- & \textbf{50} & 01:57 & -- & \textbf{O3$_{-2}^-$}, \textbf{O3$_{-1}^-$, O1} \\
    & & S3.4  & -- & \textbf{24} & 2 & 0 & \textbf{0} & \textbf{0} & -- & \textbf{50} & 04:58 & -- & \textbf{O3$_{-2}^-$}, \textbf{O3$_{-1}^-$}, \textbf{O3$_{2}^+$}, o3$_{1}^+$ \\
    & & S4.1* & -- & \textbf{24} & 2 & 0 & \textbf{0} & \textbf{0} & -- & \textbf{50} & 00:02 & -- & \textbf{O1} \\
    & & S4.2  & -- & \textbf{24} & 2 & 0 & \textbf{0} & \textbf{0} & -- & \textbf{50} & 03:53 & -- & \textbf{O3$_{-2}^-$}, \textbf{O3$_{2}^+$}, \textbf{O3$_{-1}^-$}, o3$_1^+$ \\

\pagebreak

\multirow{22}{*}{\rotatebox{90}{\textbf{D3}}}
    & \multicolumn{2}{l|@{\xhspace}}{\textit{total}} & -- & \textit{26} & \textit{49} & \textit{899} & \textit{13} & \textit{5} & -- & & & & \\
    & \multicolumn{2}{l|@{\xhspace}}{manual} & -- & 14 & 12 & 6 & \textbf{0} & \textbf{0} & -- & 40 & & & \\ 
\cmidrule(r){2-14}
& \multirow{10}{*}{\rotatebox{90}{1h}}
      & S1.1  & -- & 1 & 2 & 29 & \textbf{0} & \textbf{0} & -- & 4 & 00:23 & 00:03 & \textbf{O2} \\
    & & S1.2  & -- & 0 & 3 & 29 & \textbf{0} & \textbf{0} & -- & 3 & 00:01 & 00:01 & \textbf{O3$_{-2}^-$}, \textbf{O3$_{-1}^-$} \\
    & & S2.1  & -- & \textbf{16} & 9 & 7 & \textbf{0} & \textbf{0} & -- & \textbf{41} & 00:01 & 00:01 & \textbf{O1} \\
    & & S2.2  & -- & \textbf{16} & 9 & 7 & \textbf{0} & \textbf{0} & -- & \textbf{41} & 00:03 & 00:01 & \textbf{O3$_{2}^+$}, \textbf{O3$_{1}^+$} \\
    & & S3.1  & -- & \textbf{16} & 9 & 7 & \textbf{0} & \textbf{0} & -- & \textbf{41} & 01:59 & 00:41 & \textbf{O2}, \textbf{O1} \\
    & & S3.2  & -- & \textbf{16} & 9 & 7 & \textbf{0} & \textbf{0} & -- & \textbf{41} & 01:02 & 00:41 & \textbf{O2}, \textbf{O3$_{2}^+$}, \textbf{O3$_{1}^+$} \\
    & & S3.3  & -- & 15 & 11 & 6 & \textbf{0} & \textbf{0} & -- & \textbf{41} & 00:06 & 00:02 & \textbf{O3$_{-2}^-$}, \textbf{O3$_{-1}^-$, O1} \\
    & & S3.4  & -- & \textbf{16} & 9 & 7 & \textbf{0} & \textbf{0} & -- & \textbf{41} & 00:02 & 00:01 & \textbf{O3$_{-2}^-$}, \textbf{O3$_{-1}^-$}, \textbf{O3$_{2}^+$}, \textbf{O3$_{1}^+$} \\
    & & S4.1* & -- & \textbf{16} & 9 & 7 & \textbf{0} & \textbf{0} & -- & \textbf{41} & 00:01 & 00:01 & \textbf{O1} \\
    & & S4.2  & -- & \textbf{16} & 9 & 7 & \textbf{0} & \textbf{0} & -- & \textbf{41} & 00:03 & 00:01 & \textbf{O3$_{-2}^-$}, \textbf{O3$_{2}^+$}, \textbf{O3$_{-1}^-$}, \textbf{O3$_{1}^+$} \\
\cmidrule(r){2-14}
& \multirow{10}{*}{\rotatebox{90}{15min}}
      & S1.1  & -- & 1 & 2 & 29 & \textbf{0} & \textbf{0} & -- & 4 & 00:23 & -- & \textbf{O2} \\
    & & S1.2  & -- & 0 & 3 & 29 & \textbf{0} & \textbf{0} & -- & 3 & 00:01 & -- & \textbf{O3$_{-2}^-$}, \textbf{O3$_{-1}^-$} \\
    & & S2.1  & -- & \textbf{16} & 9 & 7 & \textbf{0} & \textbf{0} & -- & \textbf{41} & 00:01 & -- & \textbf{O1} \\
    & & S2.2  & -- & \textbf{16} & 9 & 7 & \textbf{0} & \textbf{0} & -- & \textbf{41} & 00:03 & -- & \textbf{O3$_{2}^+$}, \textbf{O3$_{1}^+$} \\
    & & S3.1  & -- & \textbf{16} & 9 & 7 & \textbf{0} & \textbf{0} & -- & \textbf{41} & 01:59 & -- & \textbf{O2}, \textbf{O1} \\
    & & S3.2  & -- & \textbf{16} & 9 & 7 & \textbf{0} & \textbf{0} & -- & \textbf{41} & 01:02 & -- & \textbf{O2}, \textbf{O3$_{2}^+$}, \textbf{O3$_{1}^+$} \\
    & & S3.3  & -- & 15 & 11 & 6 & \textbf{0} & \textbf{0} & -- & \textbf{41} & 00:06 & -- & \textbf{O3$_{-2}^-$}, \textbf{O3$_{-1}^-$, O1} \\
    & & S3.4  & -- & \textbf{16} & 9 & 7 & \textbf{0} & \textbf{0} & -- & \textbf{41} & 00:02 & -- & \textbf{O3$_{-2}^-$}, \textbf{O3$_{-1}^-$}, \textbf{O3$_{2}^+$}, \textbf{O3$_{1}^+$} \\
    & & S4.1* & -- & \textbf{16} & 9 & 7 & \textbf{0} & \textbf{0} & -- & \textbf{41} & 00:01 & -- & \textbf{O1} \\
    & & S4.2  & -- & \textbf{16} & 9 & 7 & \textbf{0} & \textbf{0} & -- & \textbf{41} & 00:03 & -- & \textbf{O3$_{-2}^-$}, \textbf{O3$_{2}^+$}, \textbf{O3$_{-1}^-$}, \textbf{O3$_{1}^+$} \\
\midrule

\multirow{22}{*}{\rotatebox{90}{\textbf{D4}}}
    & \multicolumn{2}{l|@{\xhspace}}{\textit{total}} & \textit{37} & \textit{207} & \textit{253} & \textit{55} & \textit{319} & \textit{121} & -- & & & & \\
    & \multicolumn{2}{l|@{\xhspace}}{manual*} & 27 & 26 & 18 & 4 & 32 & 9 & -- & 128 & & & \\
\cmidrule(r){2-14}
& \multirow{10}{*}{\rotatebox{90}{1h}}
      & S1.1 & 4 & 29 & 38 & 6 & 39 & \textbf{0} & -- & 73 & 00:23 & 00:05 & \textbf{O2} \\
    & & S1.2 & 6 & 35 & \textbf{56} & 12 & 7 & \textbf{0} & -- & 143 & 00:01 & 00:01 & \textbf{O3$_{-2}^-$}, \textbf{O3$_{-1}^-$} \\
    & & S2.1 & 21 & \textbf{58} & 23 & 4 & \textbf{4} & 6 & -- & \textbf{207} & 00:48 & 00:38 & \textbf{O1} \\
    & & S2.2 & \textbf{32} & 22 & 18 & 3 & 30 & 11 & -- & 138 & 00:17 & 00:16 & \textbf{O3$_{4}^+$}, \textbf{O3$_{2}^+$}, \textbf{O3$_{1}^+$} \\
    & & S3.1 & 19 & 52 & 29 & 6 & 10 & \textbf{0} & -- & 199 & 60:01 & 12:17 & \textbf{O2}, o1 \\
    & & S3.2 & 28 & 30 & 18 & 4 & 36 & \textbf{0} & -- & 154 & 01:00 & 00:49 & \textbf{O2}, \textbf{O3$_{4}^+$}, \textbf{O3$_{2}^+$}, \textbf{O3$_{1}^+$} \\
    & & S3.3 & 15 & 53 & 33 & 8 & 7 & \textbf{0} & -- & 192 & 00:32 & 00:28 & \textbf{O3$_{-2}^-$}, \textbf{O3$_{-1}^-$}, \textbf{O1} \\
    & & S3.4 & 16 & 41 & 46 & 6 & 7 & \textbf{0} & -- & 185 & 00:46 & 00:22 & \textbf{O3$_{-2}^-$}, \textbf{O3$_{-1}^-$}, \textbf{O3$_{4}^+$}, \textbf{O3$_{2}^+$}, \textbf{O3$_{1}^+$} \\
    & & S4.1 & \textbf{32} & 22 & 18 & 3 & 34 & 7 & -- & 142 & 00:12 & 00:08 & \textbf{O3$_{4}^+$}, \textbf{O1} \\
    & & S4.2 & \textbf{32} & 22 & 18 & 3 & 34 & 7 & -- & 142 & 00:35 & 00:27 & \textbf{O3$_{4}^+$}, \textbf{O3$_{-2}^-$}, \textbf{O3$_{2}^+$}, \textbf{O3$_{-1}^-$}, \textbf{O3$_{1}^+$} \\
\cmidrule(r){2-14}
& \multirow{10}{*}{\rotatebox{90}{15min}}
      & S1.1 & 4 & 29 & 38 & 6 & 39 & \textbf{0} & -- & 73 & 00:23 & -- & \textbf{O2} \\
    & & S1.2 & 6 & 35 & \textbf{56} & 12 & 7 & \textbf{0} & -- & 143 & 00:01 & -- & \textbf{O3$_{-2}^-$}, \textbf{O3$_{-1}^-$} \\
    & & S2.1 & 21 & 58 & 23 & 4 & \textbf{4} & 6 & -- & \textbf{207} & 00:48 & -- & \textbf{O1} \\
    & & S2.2 & \textbf{32} & 22 & 18 & 3 & 30 & 11 & -- & 138 & 00:17 & -- & \textbf{O3$_{4}^+$}, \textbf{O3$_{2}^+$}, \textbf{O3$_{1}^+$} \\
    & & S3.1 & 15 & \textbf{60} & 27 & 4 & 10 & \textbf{0} & -- & 197 & 07:54 & -- & \textbf{O2}, o1 \\
    & & S3.2 & 28 & 30 & 18 & 4 & 36 & \textbf{0} & -- & 154 & 01:00 & -- & \textbf{O2}, \textbf{O3$_{4}^+$}, \textbf{O3$_{2}^+$}, \textbf{O3$_{1}^+$} \\
    & & S3.3 & 15 & 53 & 33 & 8 & 7 & \textbf{0} & -- & 192 & 00:32 & -- & \textbf{O3$_{-2}^-$}, \textbf{O3$_{-1}^-$}, \textbf{O1} \\
    & & S3.4 & 16 & 41 & 46 & 6 & 7 & \textbf{0} & -- & 185 & 00:46 & -- & \textbf{O3$_{-2}^-$}, \textbf{O3$_{-1}^-$}, \textbf{O3$_{4}^+$}, \textbf{O3$_{2}^+$}, \textbf{O3$_{1}^+$} \\
    & & S4.1 & \textbf{32} & 22 & 18 & 3 & 34 & 7 & -- & 142 & 00:12 & -- & \textbf{O3$_{4}^+$}, \textbf{O1} \\
    & & S4.2 & \textbf{32} & 22 & 18 & 3 & 34 & 7 & -- & 142 & 00:35 & -- & \textbf{O3$_{4}^+$}, \textbf{O3$_{-2}^-$}, \textbf{O3$_{2}^+$}, \textbf{O3$_{-1}^-$}, \textbf{O3$_{1}^+$} \\

\pagebreak

\multirow{22}{*}{\rotatebox{90}{\textbf{D5}}}
    & \multicolumn{2}{l|@{\xhspace}}{\textit{total}} & \textit{6} & \textit{110} & \textit{285} & \textit{525} & \textit{344} & \textit{52} & \textit{10} & & & & \\
    & \multicolumn{2}{l|@{\xhspace}}{manual*} & 5 & 14 & 17 & 46 & 18 & 2 & \textbf{0} & 43 & & & \\
\cmidrule(r){2-14}
& \multirow{10}{*}{\rotatebox{90}{1h}}
      & S1.1 & 0 & 14 & 30 & 58 & \textbf{0} & \textbf{0} & \textbf{0} & 58 & 00:34 & 00:08 & \textbf{O2} \\
    & & S1.2 & 1 & 14 & 33 & 54 & \textbf{0} & \textbf{0} & \textbf{0} & 65 & 00:01 & 00:01 & \textbf{O3$_{-4}^-$}, \textbf{O3$_{-2}^-$}, \textbf{O3$_{-1}^-$} \\
    & & S2.1 & 3 & \textbf{50} & 39 & 8 & 2 & \textbf{0} & \textbf{0} &\textbf{149} & 00:18 & 00:06 & \textbf{O1} \\
    & & S2.2 & \textbf{6} & 46 & 35 & 11 & 4 & \textbf{0} & \textbf{0} & 147 & 60:01 & 00:17 & \textbf{O3$_{4}^+$}, \textbf{O3$_{2}^+$}, o3$_{1}^+$ \\
    & & S3.1 & 3 & 44 & \textbf{47} & 8 & \textbf{0} & \textbf{0} & \textbf{0} & 147 & 60:01 & 01:31 & \textbf{O2}, o1 \\
    & & S3.2 & \textbf{6} & 46 & 23 & 27 & \textbf{0} & \textbf{0} & \textbf{0} & 139 & 60:01 & 01:07 & \textbf{O2}, \textbf{O3$_{4}^+$}, o3$_{2}^+$, \textit{o3$_{1}^+$} \\
    & & S3.3 & 3 & 44 & \textbf{47} & 8 & \textbf{0} & \textbf{0} & \textbf{0} & 147 & 60:01 & 01:09 & \textbf{O3$_{-4}^-$}, \textbf{O3$_{-2}^-$}, \textbf{O3$_{-1}^-$}, o1 \\
    & & S3.4 & \textbf{6} & 46 & 29 & 21 & \textbf{0} & \textbf{0} & \textbf{0} & 145 & 60:01 & 00:31 & \textbf{O3$_{-4}^-$}, \textbf{O3$_{-2}^-$}, \textbf{O3$_{-1}^-$}, \textbf{O3$_{4}^+$}, \textbf{O3$_{2}^+$}, o3$_{1}^+$ \\
    & & S4.1 & \textbf{6} & 44 & 37 & 13 & 2 & \textbf{0} & \textbf{0} & 147 & 00:16 & 00:06 & \textbf{O3$_{-4}^-$}, \textbf{O3$_{4}^+$}, \textbf{O1} \\
    & & S4.2 & \textbf{6} & 46 & 29 & 21 & \textbf{0} & \textbf{0} & \textbf{0} & 145 & 60:01 & 00:20 & \textbf{O3$_{-4}^-$}, \textbf{O3$_{4}^+$}, \textbf{O3$_{-2}^-$}, \textbf{O3$_{2}^+$}, \textbf{O3$_{-1}^-$}, o3$_{1}^+$ \\
\cmidrule(r){2-14}
& \multirow{10}{*}{\rotatebox{90}{15min}}
      & S1.1 & 0 & 14 & 30 & 58 & \textbf{0} & \textbf{0} & \textbf{0} & 58 & 00:34 & -- & \textbf{O2} \\
    & & S1.2 & 1 & 14 & 33 & 54 & \textbf{0} & \textbf{0} & \textbf{0} & 65 & 00:01 & -- & \textbf{O3$_{-4}^-$}, \textbf{O3$_{-2}^-$}, \textbf{O3$_{-1}^-$} \\
    & & S2.1 & 3 & \textbf{50} & 39 & 8 & 2 & \textbf{0} & \textbf{0} & \textbf{149} & 00:18 & -- & \textbf{O1} \\
    & & S2.2 & \textbf{6} & 46 & 35 & 11 & 4 & \textbf{0} & \textbf{0} & 147 & 05:05 & -- & \textbf{O3$_{4}^+$}, \textbf{O3$_{2}^+$}, o3$_{1}^+$ \\
    & & S3.1 & 3 & 44 & \textbf{47} & 8 & \textbf{0} & \textbf{0} & \textbf{0} & 147 & 08:14 & -- & \textbf{O2}, o1 \\
    & & S3.2 & \textbf{6} & 46 & 29 & 21 & \textbf{0} & \textbf{0} & \textbf{0} & 145 & 08:28 & -- & \textbf{O2}, \textbf{O3$_{4}^+$}, o3$_{2}^+$, o3$_{1}^+$ \\
    & & S3.3 & 3 & 44 & \textbf{47} & 8 & \textbf{0} & \textbf{0} & \textbf{0} & 147 & 03:48 & -- & \textbf{O3$_{-4}^-$}, \textbf{O3$_{-2}^-$}, \textbf{O3$_{-1}^-$}, o1 \\
    & & S3.4 & \textbf{6} & 46 & 29 & 21 & \textbf{0} & \textbf{0} & \textbf{0} & 145 & 02:56 & -- & \textbf{O3$_{-4}^-$}, \textbf{O3$_{-2}^-$}, \textbf{O3$_{-1}^-$}, \textbf{O3$_{4}^+$}, \textbf{O3$_{2}^+$}, o3$_{1}^+$ \\
    & & S4.1 & \textbf{6} & 44 & 37 & 13 & 2 & \textbf{0} & \textbf{0} & 147 & 00:16 & -- & \textbf{O3$_{-4}^-$}, \textbf{O3$_{4}^+$}, \textbf{O1} \\
    & & S4.2 & \textbf{6} & 46 & 29 & 21 & \textbf{0} & \textbf{0} & \textbf{0} & 145 & 02:35 & -- & \textbf{O3$_{-4}^-$}, \textbf{O3$_{4}^+$}, \textbf{O3$_{-2}^-$}, \textbf{O3$_{2}^+$}, \textbf{O3$_{-1}^-$}, o3$_{1}^+$ \\
\midrule

\multirow{22}{*}{\rotatebox{90}{\textbf{D6}}}
    & \multicolumn{2}{l|@{\xhspace}}{\textit{total}} & \textit{138} & \textit{96} & \textit{81} & \textit{3018} & \textit{94} & \textit{113} & -- & & & & \\
    & \multicolumn{2}{l|@{\xhspace}}{manual} & 77 & 14 & 6 & 134 & 5 & 4 & -- & 329 & & & \\
\cmidrule(r){2-14}
& \multirow{10}{*}{\rotatebox{90}{1h}}
      & S1.1 & 4 & 4 & 3 & 237 & \textbf{0} & \textbf{0} & -- & 27 & 01:28 & 00:39 & \textbf{O2} \\
    & & S1.2 & 12 & 5 & 6 & 223 & \textbf{0} & \textbf{0} & -- & 64 & 00:01 & 00:01 & \textbf{O3$_{-2}^-$}, \textbf{O3$_{-1}^-$} \\
    & & S2.1 & 106 & 28 & 8 & 103 & 2 & 1 & -- & \textbf{484} & 04:12 & 01:03 & \textbf{O1} \\
    & & S2.2 & \textbf{107} & 26 & 7 & 105 & 2 & 3 & -- & 479 & 12:56 & 03:53 & \textbf{O3$_{4}^+$}, \textbf{O3$_{2}^+$}, \textbf{O3$_{1}^+$} \\
    & & S3.1 & 102 & 27 & 8 & 111 & \textbf{0} & \textbf{0} & -- & 470 & 34:10 & 11:54 & \textbf{O2}, \textbf{O1} \\
    & & S3.2 & 103 & 24 & 6 & 115 & \textbf{0} & \textbf{0} & -- & 466 & 21:30 & 13:33 & \textbf{O2}, \textbf{O3$_{4}^+$}, \textbf{O3$_{2}^+$}, \textbf{O3$_{1}^+$} \\
    & & S3.3 & 99 & \textbf{32} & \textbf{10} & 109 & \textbf{0} & \textbf{0} & -- & 470 & 00:43 & 00:27 & \textbf{O3$_{-2}^-$}, \textbf{O3$_{-1}^-$}, \textbf{O1} \\
    & & S3.4 & 103 & 24 & 6 & 117 & \textbf{0} & \textbf{0} & -- & 466 & 16:40 & 01:09 & \textbf{O3$_{-2}^-$}, \textbf{O3$_{-1}^-$}, \textbf{O3$_{4}^+$}, \textbf{O3$_{2}^+$}, \textbf{O3$_{1}^+$} \\
    & & S4.1 & \textbf{107} & 25 & 6 & 107 & 2 & 1 & -- & 480 & 04:26 & 02:18 & \textbf{O3$_{4}^+$}, \textbf{O1} \\
    & & S4.2 & \textbf{107} & 25 & 6 & 107 & 2 & 1 & -- & 480 & 04:30 & 02:47 & \textbf{O3$_{4}^+$}, \textbf{O3$_{-2}^-$}, \textbf{O3$_{2}^+$}, \textbf{O3$_{-1}^-$}, \textbf{O3$_{1}^+$} \\
\cmidrule(r){2-14}
& \multirow{10}{*}{\rotatebox{90}{15min}}
      & S1.1 & 4 & 4 & 3 & 237 & \textbf{0} & \textbf{0} & -- & 27 & 01:28 & -- & \textbf{O2} \\
    & & S1.2 & 12 & 5 & 6 & 223 & \textbf{0} & \textbf{0} & -- & 64 & 00:01 & -- & \textbf{O3$_{-2}^-$}, \textbf{O3$_{-1}^-$} \\
    & & S2.1 & 106 & 28 & 8 & 103 & 2 & 1 & -- & \textbf{484} & 04:12 & -- & \textbf{O1} \\
    & & S2.2 & \textbf{107} & 26 & 7 & 105 & 2 & 3 & -- & 479 & 08:49 & -- & \textbf{O3$_{4}^+$}, \textbf{O3$_{2}^+$}, o3$_{1}^+$ \\
    & & S3.1 & 102 & 27 & 8 & 111 & \textbf{0} & \textbf{0} & -- & 470 & 09:15 & -- & \textbf{O2}, o1 \\
    & & S3.2 & 103 & 24 & 6 & 115 & \textbf{0} & \textbf{0} & -- & 466 & 12:43 & -- & \textbf{O2}, o3$_{4}^+$, o3$_{2}^+$, o3$_{1}^+$ \\
    & & S3.3 & 99 & \textbf{32} & \textbf{10} & 109 & \textbf{0} & \textbf{0} & -- & 470 & 00:43 & -- & \textbf{O3$_{-2}^-$}, \textbf{O3$_{-1}^-$}, \textbf{O1} \\
    & & S3.4 & 103 & 24 & 6 & 117 & \textbf{0} & \textbf{0} & -- & 466 & 06:33 & -- & \textbf{O3$_{-2}^-$}, \textbf{O3$_{-1}^-$}, \textbf{O3$_{4}^+$}, o3$_{2}^+$, o3$_{1}^+$ \\
    & & S4.1 & \textbf{107} & 25 & 6 & 107 & 2 & 1 & -- & 480 & 04:26 & -- & \textbf{O3$_{4}^+$}, \textbf{O1} \\
    & & S4.2 & \textbf{107} & 25 & 6 & 107 & 2 & 1 & -- & 480 & 04:30 & -- & \textbf{O3$_{4}^+$}, \textbf{O3$_{-2}^-$}, \textbf{O3$_{2}^+$}, \textbf{O3$_{-1}^-$}, \textbf{O3$_{1}^+$} \\

\pagebreak

\multirow{22}{*}{\rotatebox{90}{\textbf{D7}}}
    & \multicolumn{2}{l|@{\xhspace}}{\textit{total}} & \textit{4} & \textit{43} & \textit{1490} & \textit{2272} & \textit{881} & \textit{2} & -- & & & & \\
    & \multicolumn{2}{l|@{\xhspace}}{manual*} & 4 & 22 & 61 & 82 & 28 & 1 & -- & 91 & & & \\
\cmidrule(r){2-14}
& \multirow{10}{*}{\rotatebox{90}{1h}}
      & S1.1 & 2 & 1 & 88 & 107 & \textbf{0} & \textbf{0} & -- & 98 & 06:10 & 01:14 & \textbf{O2} \\
    & & S1.2 & 0 & 0 & 80 & 118 & \textbf{0} & \textbf{0} & -- & 80 & 00:01 & 00:01 & \textbf{O3$_{-2}^-$}, \textbf{O3$_{-1}^-$} \\
    & & S2.1 & \textbf{4} & 14 & 173 & 7 & \textbf{0} & \textbf{0} & -- & \textbf{217} & 60:01 & 52:47 & o1 \\
    & & S2.2 & \textbf{4} & \textbf{31} & 119 & 26 & 17 & 1 & -- & 178 & 60:01 & 33:11 & \textbf{O3$_{4}^+$}, \textbf{O3$_{2}^+$}, o3$_{1}^+$ \\
    & & S3.1 & 2 & 5 & 186 & 5 & \textbf{0} & \textbf{0} & -- & 204 & 60:01 & 50:09 & \textbf{O2}, o1 \\
    & & S3.2 & \textbf{4} & 27 & 53 & 114 & \textbf{0} & \textbf{0} & -- & 123 & 60:01 & 43:01 & \textbf{O2}, \textbf{O3$_{4}^+$}, o3$_{2}^+$, \textit{o3$_{1}^+$} \\
    & & S3.3 & \textbf{4} & 10 & \textbf{177} & 7 & \textbf{0} & \textbf{0} & -- & 213 & 60:01 & 57:21 & \textbf{O3$_{-2}^-$}, \textbf{O3$_{-1}^-$}, o1 \\
    & & S3.4 & \textbf{4} & 27 & 135 & 32 & \textbf{0} & \textbf{0} & -- & 205 & 60:01 & 37:33 & \textbf{O3$_{-2}^-$}, \textbf{O3$_{-1}^-$}, \textbf{O3$_{4}^+$}, \textbf{O3$_{2}^+$}, o3$_{1}^+$ \\
    & & S4.1 & \textbf{4} & 14 & 171 & 9 & \textbf{0} & \textbf{0} & -- & 215 & 60:01 & 42:01 & \textbf{O3$_{4}^+$}, o1 \\
    & & S4.2 & \textbf{4} & 30 & 119 & 40 & 5 & \textbf{0} & -- & 190 & 60:01 & 46:42 & \textbf{O3$_{4}^+$}, \textbf{O3$_{-2}^-$}, \textbf{O3$_{2}^+$}, \textbf{O3$_{-1}^-$}, o3$_{1}^+$ \\
\cmidrule(r){2-14}
& \multirow{10}{*}{\rotatebox{90}{15min}}
      & S1.1 & 2 & 1 & 88 & 107 & \textbf{0} & \textbf{0} & -- & 98 & 06:10 & -- & \textbf{O2} \\
    & & S1.2 & 0 & 0 & 80 & 118 & \textbf{0} & \textbf{0} & -- & 80 & 00:01 & -- & \textbf{O3$_{-2}^-$}, \textbf{O3$_{-1}^-$} \\
    & & S2.1 & 2 & 19 & 157 & 20 & \textbf{0} & \textbf{0} & -- & \textbf{203} & 15:01 & -- & o1 \\
    & & S2.2 & \textbf{4} & \textbf{31} & 103 & 36 & 23 & 1 & -- & 156 & 05:05 & -- & \textbf{O3$_{4}^+$}, \textbf{O3$_{2}^+$}, o3$_{1}^+$ \\
    & & S3.1 & \textbf{4} & 2 & 137 & 55 & \textbf{0} & \textbf{0} & -- & 157 & 14:46 & -- & \textbf{O2}, o1 \\
    & & S3.2 & \textbf{4} & 2 & 142 & 50 & \textbf{0} & \textbf{0} & -- & 162 & 12:26 & -- & o2, \textbf{O3$_{4}^+$}, o3$_{2}^+$, o3$_{1}^+$ \\
    & & S3.3 & \textbf{4} & 4 & \textbf{171} & 19 & \textbf{0} & \textbf{0} & -- & 195 & 05:15 & -- & \textbf{O3$_{-2}^-$}, \textbf{O3$_{-1}^-$}, o1 \\
    & & S3.4 & \textbf{4} & 27 & 119 & 48 & \textbf{0} & \textbf{0} & -- & 189 & 06:04 & -- & \textbf{O3$_{-2}^-$}, \textbf{O3$_{-1}^-$}, \textbf{O3$_{4}^+$}, o3$_{2}^+$, o3$_{1}^+$ \\
    & & S4.1 & \textbf{4} & 8 & 132 & 54 & \textbf{0} & \textbf{0} & -- & 164 & 07:39 & -- & \textbf{O3$_{4}^+$}, o1 \\
    & & S4.2 & \textbf{4} & 30 & 87 & 72 & 5 & \textbf{0} & -- & 158 & 03:15 & -- & \textbf{O3$_{4}^+$}, \textbf{O3$_{-2}^-$}, \textbf{O3$_{2}^+$}, \textbf{O3$_{-1}^-$}, o3$_{1}^+$ \\
\midrule

\multirow{22}{*}{\rotatebox{90}{\textbf{D8.1}}}
    & \multicolumn{2}{l|@{\xhspace}}{\textit{total}} & \textit{23} & \textit{25} & \textit{1923} & \textit{2961} & \textit{1223} & \textit{7} & -- & & & & \\
    & \multicolumn{2}{l|@{\xhspace}}{manual*} & 17 & 11 & 46 & 55 & 23 & 2 & -- & 109 & & & \\
\cmidrule(r){2-14}
& \multirow{10}{*}{\rotatebox{90}{1h}}
      & S1.1 & 6 & 0 & \textbf{148} & 0 & \textbf{0} & \textbf{0} & -- & 172 & 25:59 & 11:44 & \textbf{O2} \\
    & & S1.2 & 0 & 0 & 55 & 99 & \textbf{0} & \textbf{0} & -- & 55 & 00:01 & 00:01 & \textbf{O3$_{-2}^-$}, \textbf{O3$_{-1}^-$} \\
    & & S2.1 & \textbf{23} & 7 & 114 & 9 & 1 & \textbf{0} & -- & \textbf{219} & 60:01 & 37:08 & o1 \\
    & & S2.2 & \textbf{23} & \textbf{15} & 95 & 10 & 10 & 1 & -- & 205 & 60:01 & 22:37 & \textbf{O3$_{4}^+$}, \textbf{O3$_{2}^+$}, o3$_{1}^+$ \\
    & & S3.1 & 6 & 1 & 147 & 0 & \textbf{0} & \textbf{0} & -- & 173 & 60:01 & 41:22 & \textbf{O2}, o1 \\
    & & S3.2 & 7 & 0 & 147 & 0 & \textbf{0} & \textbf{0} & -- & 175 & 60:01 & 38:15 & \textbf{O2}, o3$_4^+$, \textit{o3$_{2}^+$, o3$_{1}^+$} \\ 
    & & S3.3 & 22 & 5 & 120 & 7 & \textbf{0} & \textbf{0} & -- & 218 & 60:01 & 47:12 & \textbf{O3$_{-2}^-$}, \textbf{O3$_{-1}^-$}, o1 \\
    & & S3.4 & 22 & 12 & 100 & 20 & \textbf{0} & \textbf{0} & -- & 212 & 60:01 & 38:17 & \textbf{O3$_{-2}^-$}, \textbf{O3$_{-1}^-$}, \textbf{O3$_{4}^+$}, \textbf{O3$_{2}^+$}, o3$_{1}^+$ \\
    & & S4.1 & \textbf{23} & 1 & 112 & 17 & 1 & \textbf{0} & -- & 205 & 60:01 & 57:23 & \textbf{O3$_{4}^+$}, o1 \\
    & & S4.2 & \textbf{23} & 14 & 96 & 18 & 3 & \textbf{0} & -- & 213 & 60:01 & 54:22 & \textbf{O3$_{4}^+$}, \textbf{O3$_{-2}^-$}, \textbf{O3$_{2}^+$}, \textbf{O3$_{-1}^-$}, o3$_{1}^+$ \\
\cmidrule(r){2-14}
& \multirow{10}{*}{\rotatebox{90}{15min}}
      & S1.1 & 6 & 0 & \textbf{148} & 0 & \textbf{0} & \textbf{0} & -- & 172 & 15:01 & -- & o2 \\
    & & S1.2 & 0 & 0 & 55 & 99 & \textbf{0} & \textbf{0} & -- & 55 & 00:01 & -- & \textbf{O3$_{-2}^-$}, \textbf{O3$_{-1}^-$} \\
    & & S2.1 & 22 & 8 & 105 & 19 & \textbf{0} & \textbf{0} & -- & \textbf{209} & 15:01 & -- & o1 \\
    & & S2.2 & \textbf{23} & \textbf{15} & 77 & 29 & 9 & 1 & -- & 188 & 05:08 & -- & \textbf{O3$_{4}^+$}, \textbf{O3$_{2}^+$}, o3$_{1}^+$ \\
    & & S3.1 & 4 & 1 & 50 & 75 & 23 & 1 & -- & 43 & 19:03 & -- & o2, o1 \\
    & & S3.2 & 4 & 2 & 98 & 34 & 16 & \textbf{0} & -- & 102 & 16:09 & -- & o2, o3$_4^+$, o3$_{2}^+$, o3$_{1}^+$ \\
    & & S3.3 & 20 & 1 & 121 & 12 & \textbf{0} & \textbf{0} & -- & 203 & 05:39 & -- & \textbf{O3$_{-2}^-$}, \textbf{O3$_{-1}^-$}, o1 \\
    & & S3.4 & 22 & 12 & 83 & 37 & \textbf{0} & \textbf{0} & -- & 195 & 03:47 & -- & \textbf{O3$_{-2}^-$}, \textbf{O3$_{-1}^-$}, \textbf{O3$_{4}^+$}, \textbf{O3$_{2}^+$}, o3$_{1}^+$ \\
    & & S4.1 & \textbf{23} & 1 & 102 & 25 & 3 & \textbf{0} & -- & 193 & 07:49 & -- & \textbf{O3$_{4}^+$}, o1 \\
    & & S4.2 & \textbf{23} & 14 & 78 & 36 & 3 & \textbf{0} & -- & 195 & 03:51 & -- & \textbf{O3$_{4}^+$}, \textbf{O3$_{-2}^-$}, \textbf{O3$_{2}^+$}, \textbf{O3$_{-1}^-$}, o3$_{1}^+$ \\

\pagebreak

\multirow{22}{*}{\rotatebox{90}{\textbf{D8.2}}}
    & \multicolumn{2}{l|@{\xhspace}}{\textit{total}} & \textit{23} & \textit{25} & \textit{50} & \textit{6036} & \textit{21} & \textit{7} & -- & & & & \\
    & \multicolumn{2}{l|@{\xhspace}}{manual} & 17 & 11 & 10 & 111 & 3 & 2 & -- & 93 & & & \\
\cmidrule(r){2-14}
& \multirow{10}{*}{\rotatebox{90}{1h}}
      & S1.1 & 0 & 1 & 2 & 151 & \textbf{0} & \textbf{0} & -- & 4 & 06:54 & 01:14 & \textbf{O2} \\
    & & S1.2 & 0 & 3 & 1 & 150 & \textbf{0} & \textbf{0} & -- & 7 & 00:01 & 00:01 & \textbf{O3$_{-2}^-$}, \textbf{O3$_{-1}^-$} \\
    & & S2.1 & \textbf{23} & 14 & 16 & 101 & \textbf{0} & \textbf{0} & -- & \textbf{136} & 00:25 & 00:15 & \textbf{O1} \\
    & & S2.2 & \textbf{23} & \textbf{15} & \textbf{17} & 96 & 1 & 2 & -- & 134 & 00:31 & 00:11 & \textbf{O3$_{4}^+$}, \textbf{O3$_{2}^+$}, \textbf{O3$_{1}^+$} \\
    & & S3.1 & 2 & 2 & 0 & 150 & \textbf{0} & \textbf{0} & -- & 12 & 60:01 & 08:17 & \textbf{O2}, o1 \\
    & & S3.2 & 4 & 0 & 0 & 150 & \textbf{0} & \textbf{0} & -- & 16 & 60:01 & 11:51 & \textbf{O2}, o3$_{4}^+$, \textit{o3$_{2}^+$}, \textit{o3$_{1}^+$} \\
    & & S3.3 & \textbf{23} & 14 & 16 & 101 & \textbf{0} & \textbf{0} & -- & \textbf{136} & 01:17 & 00:45 & \textbf{O3$_{-2}^-$}, \textbf{O3$_{-1}^-$}, \textbf{O1} \\
    & & S3.4 & \textbf{23} & 14 & 16 & 101 & \textbf{0} & \textbf{0} & -- & \textbf{136} & 00:19 & 00:14 & \textbf{O3$_{-2}^-$}, \textbf{O3$_{-1}^-$}, \textbf{O3$_{4}^+$}, \textbf{O3$_{2}^+$}, \textbf{O3$_{1}^+$} \\
    & & S4.1 & \textbf{23} & 14 & 16 & 101 & \textbf{0} & \textbf{0} & -- & \textbf{136} & 00:46 & 00:37 & \textbf{O3$_{4}^+$}, \textbf{O1} \\
    & & S4.2 & \textbf{23} & 14 & 16 & 101 & \textbf{0} & \textbf{0} & -- & \textbf{136} & 00:18 & 00:08 & \textbf{O3$_{4}^+$}, \textbf{O3$_{-2}^-$}, \textbf{O3$_{2}^+$}, \textbf{O3$_{-1}^-$}, \textbf{O3$_{1}^+$} \\
\cmidrule(r){2-14}
& \multirow{10}{*}{\rotatebox{90}{15min}}
      & S1.1 & 0 & 1 & 2 & 151 & \textbf{0} & \textbf{0} & -- & 4 & 06:54 & -- & \textbf{O2} \\
    & & S1.2 & 0 & 3 & 1 & 150 & \textbf{0} & \textbf{0} & -- & 7 & 00:01 & -- & \textbf{O3$_{-2}^-$}, \textbf{O3$_{-1}^-$} \\
    & & S2.1 & \textbf{23} & 14 & 16 & 101 & \textbf{0} & \textbf{0} & -- & \textbf{136} & 00:25 & -- & \textbf{O1} \\
    & & S2.2 & \textbf{23} & \textbf{15} & \textbf{17} & 96 & 1 & 2 & -- & 134 & 00:31 & -- & \textbf{O3$_{4}^+$}, \textbf{O3$_{2}^+$}, \textbf{O3$_{1}^+$} \\
    & & S3.1 & 2 & 2 & 0 & 150 & \textbf{0} & \textbf{0} & -- & 12 & 14:44 & -- & o2, o1 \\
    & & S3.2 & 3 & 3 & 4 & 144 & \textbf{0} & \textbf{0} & -- & 22 & 15:08 & -- & o2, o3$_4^+$, o3$_{2}^+$, o3$_{1}^+$ \\
    & & S3.3 & \textbf{23} & 14 & 16 & 101 & \textbf{0} & \textbf{0} & -- & \textbf{136} & 01:17 & -- & \textbf{O3$_{-2}^-$}, \textbf{O3$_{-1}^-$}, \textbf{O1} \\
    & & S3.4 & \textbf{23} & 14 & 16 & 101 & \textbf{0} & \textbf{0} & -- & \textbf{136} & 00:19 & -- & \textbf{O3$_{-2}^-$}, \textbf{O3$_{-1}^-$}, \textbf{O3$_{4}^+$}, \textbf{O3$_{2}^+$}, \textbf{O3$_{1}^+$} \\
    & & S4.1 & \textbf{23} & 14 & 16 & 101 & \textbf{0} & \textbf{0} & -- & \textbf{136} & 00:46 & -- & \textbf{O3$_{4}^+$}, \textbf{O1} \\
    & & S4.2 & \textbf{23} & 14 & 16 & 101 & \textbf{0} & \textbf{0} & -- & \textbf{136} & 00:18 & -- & \textbf{O3$_{4}^+$}, \textbf{O3$_{-2}^-$}, \textbf{O3$_{2}^+$}, \textbf{O3$_{-1}^-$}, \textbf{O3$_{1}^+$} \\

\bottomrule

\caption{Experimental results of each strategy with a hard time limit of one hour as well as a total time limit of 15 minutes divided equally among all objectives the strategy employs.
$T$ denotes the total runtime of the computation while $T^*$ denotes the time to reach the final solution (i.e., after which point no other improvements to the solution were made by the optimizer).
The objectives are listed in the order they were optimized for with the following semantics regarding their font styles: bold upper case -- solved optimally, lower case -- optimization of this objective timed out, italic lower case -- objective was not considered during the optimization process due to a premature timeout.
We observe that in the majority of the cases, a very good solution is attained in only a fraction of the time compared to the optimal solution or the best solution attainable within the one-hour limit.}
\label{tab:exp_results}
\end{longtable}

\endgroup

\section{Discussion, Limitations and Future Work}
\label{sec:discussion}

In this section, we make a synthetic overview of our proposed approach and its outcomes, discussing a number of interesting points, limitations, and directions for future work.

\subsection{Key findings}
First, we present three key findings and overarching trends that emerged consistently across all datasets analyzed in our experiments, despite variations in course or team size, the number of skills considered, and the methods used to elicit student preferences.
These insights are particularly useful not only for researchers but also for educators seeking to implement our proposed approach in their own classrooms or educational contexts.

\subsubsection{Algorithmic solutions outperform manual ones}
Across all but two datasets, eight out of the ten strategies tested result in team formation solutions of higher quality than the manual teacher solutions.
The primary strategies producing worse solutions are S1.1 and S1.2, which focus only on avoiding the realization of negative preferences, i.e., preventing students from being placed with peers they do not wish to work with or are incompatible with.
As a result, these strategies are less effective for practical applications compared to the strategies that also incorporate positive preferences, i.e., student preferences for teammates they would like to cooperate with.
In addition to this universal observation, strategies S3.1 and S3.2 also result in subpar solutions for datasets D8.1 and D8.2.
Here, at most only the first objective O2 could be solved optimally within the time limit.
Consequently, the optimization also mainly focused on avoiding negative preferences.
This also underscores the need for more careful consideration in integrating objective O2 into the strategy design.

\subsubsection{Design strategies considering trade-offs and focusing on individual preference values to minimize dissatisfaction}
We now look at how educators can design the strategy that best fits their specific needs. 
Our results show that although many strategies produce solutions of comparable quality, no single strategy is universally the best one across all datasets.
Instead, tradeoffs need to be carefully considered. For example, one may opt to achieve the highest overall solution quality at the cost of realizing some negative preferences, or alternatively, accept slightly lower (yet still satisfactory) overall quality to ensure that no student is assigned to work with someone they prefer to avoid.
These trade-offs can also be framed as a choice between \textit{maximizing satisfaction} versus \textit{minimizing misery}: forming the best teams overall, even if some students are placed to work with undesired teammates, versus accepting a less optimal global team formation solution to prevent any such negative matching.

With this in mind, our recommendation is for teachers and students to decide carefully which metrics are most important in their context and design a suitable strategy by selecting and prioritizing the appropriate objectives accordingly.
For example, if placing students who do not want to work together into separate teams is a higher priority than teaming up those who do, a suitable strategy could start with multiple variants of objective O3 to minimize negative preference realizations, followed by objective O1 to maximize satisfaction.
This approach corresponds to our proposed strategy S3.3.

However, the space of possible strategies is not limited to those we explored in this paper.
Thanks to the hierarchical nature of our approach, new strategies can be created, e.g., strategies that prioritize maximizing satisfaction over minimizing dissatisfaction.

Further, while no single strategy proved to be the overall best one, some exhibited weaknesses that may render them unsuitable for certain contexts; especially when minimizing dissatisfaction is critical.
Our results indicate that strategies incorporating objective O2 (which explicitly minimizes teammate preference dissatisfaction) consistently required longer runtimes compared to their counterparts using O3.
If there are time constraints, the use of O2 can result in worse solutions, because the optimization may time out before all objectives have been considered, as observed in the experiments with datasets D8.1 and D8.2.
Consequently, when designing strategies that aim at minimizing dissatisfaction, we recommend avoiding objective O2 in favor of using multiple variants of objective O3 to minimize the number of unwanted realized preference values individually while remaining computationally efficient.

\subsubsection{Use timeboxing and sparse preference matrices if runtime is critical}

In practical applications, such as those encountered in educational settings, the available time is often a more important concern than achieving the absolutely optimal team formation. 

All strategies in our evaluation were given a one-hour time limit for optimization. Timeouts occurred most frequently during the optimization of the final objective, i.e., the objective that is the least important in the hierarchical order. Our results also showed that the most important, i.e., higher-priority, objectives were still solved optimally. Additionally, the most significant improvements in solution quality occurred early in the optimization process. In most cases, the solutions computed within a 15-minute time limit, equally divided across all objectives, were of similar quality to those obtained after one hour.

These findings support our recommendation that, in time-sensitive settings, it is better to opt for shorter time limits, which still yield fairly good solutions, rather than prolong the computation in search of the optimal solution. To make an even better use of the shorter time frames, we also recommend timeboxing, i.e., distributing the total available time across all objectives, to ensure that all objectives are considered in the optimization and improve the solution in cases where optimizing a single objective for too long would otherwise prevent later objectives from being considered at all. Future work could explore more nuanced time allocation approaches, such as dynamically adjusting time based on objective complexity, and examining the impact of different total time limits on the solution quality.

Another key factor affecting runtime is the sparsity of the preference matrix, especially in large courses. As shown in our experiments comparing datasets D8.1 and D8.2, using sparse preference matrices, such as those obtained from a team dating session, significantly reduces runtime compared to the use of dense matrices. In contrast, our experiments showed that for smaller courses, dense matrices can be used without adversely affecting runtime. For such cases, richer preference data, potentially derived from profile similarities, can be utilized to improve the solution quality without sacrificing runtime.
Future research could further explore the relationship between matrix sparsity and optimization performance, for instance by systematically varying sparsity levels and measuring their impact on both runtime and solution quality.

\subsection{General considerations}
We now turn to three broader remarks concerning the generalizability and contributions of our approach. 

\subsubsection{Generalization of the approach}
In designing our method, we placed intentional focus on ensuring its flexibility and adaptability across different team formation needs and settings. This generalization ability is achieved through two main algorithm design choices.

First, our method employs a modular and hierarchical approach for its different objectives.
Objective modularity is not an algorithm design choice typically adopted; most research works on algorithmic team formation optimization for a model of fixed objectives and constraints.
This modularity can be especially crucial for practical settings, where the decision-makers do not have identical team formation needs.
EDU-TF permits them to easily tailor the optimization strategy to their requirements.

Second, what we describe as ``student preferences'' can be adapted to represent different pairwise compatibility metrics, such as actual teammate preferences or profile similarities. As also seen in our paper, some educational settings allow candidate teammates to provide explicit teammate preferences, for example through a team dating process or by using prior knowledge of past collaborations. In other settings, obtaining direct teammate preferences may be impossible or undesirable due to time constraints or the size of the candidate teammate pool, such as in large courses or in MOOC settings. In these cases, pairwise preference values can be calculated from the respective profile similarities, indicating how well the specific pair can be expected to work together, based on factors relevant to the specific collaborative setting, such as personality, work style, or value alignment.

\subsubsection{Measuring team interaction quality}
Another point concerns the metric used to characterize and model the quality of team interactions.
As in most studies in the field, our model assumes that the quality of team interactions can be approximated by the sum of pairwise teammate preferences or profile similarities.
However, this approach does not fully capture the complexity of team dynamics, which may not be adequately represented by a simple sum of pairwise interactions.
Future work could therefore explore alternative metrics to assess interpersonal complementarity beyond pairwise preferences, incorporating team-level measures that consider the group as a whole.
Such an approach could offer a more nuanced understanding of team compatibility and cohesion, potentially leading to more effective team formation solutions.
From an algorithmic perspective, this could translate to exploring graph partitioning in more than two dimensions.

\subsubsection{Human agency in algorithmic team formation}
Finally, our work makes an important contribution toward increasing human agency in algorithmic team formation decision-making. 
Agency in our approach is incorporated in two ways: for students, by allowing them to explicitly express their teammate preferences, which are then integrated into the algorithmic modeling; and for teachers, by enabling them to select the specific objectives that will be used for their course's team formation optimization, alongside the ability to adjust multiple other parameters (e.g., team size, skill coverage). 
Future research could further advance this approach by allowing students and teachers to interact with the algorithm directly, for example by allowing them to adjust objectives dynamically, adjust objective weights, and control additional parameters and internal mechanics, eventually leading to a more personalized, participatory, and adaptable team formation process that better reflects user preferences and priorities.

\subsection{Limitations}
Aside from its contributions, this work also has a number of limitations. 
First, our evaluation metrics focus primarily on the number of student teammate preferences honored and computational time, rather than the impact on the academic performance of the teams.
As mentioned in section \ref{sec:results}, this evaluation choice was guided by ethical considerations present in the educational setting, where different methods for team formation could influence the students' learning outcome and grades.
Future work could nevertheless explore academic outcome measures in greater depth and aim to gather qualitative feedback from both students and teachers, providing a more holistic understanding of the algorithm’s impact on team dynamics and learning experiences.
Such insights could help further refine the algorithm and inform best practices for its application in diverse educational contexts.

A second limitation of this study is the lack of consideration for team dynamics that develop over time, such as communication styles, collaboration efficiency, and conflict resolution skills.
While our metrics primarily measure initial team preferences and setup, they do not account for how well teams function as a unit throughout the course.
Future research could explore adapting the team formation algorithms to account for these evolving dynamics, potentially by incorporating periodic feedback from team members on factors like collaboration quality and communication effectiveness into the problem modeling. 

The final limitation concerns the way that we approach the skills required to successfully complete the given project.
Our concept of skill coverage per team is relatively basic, requiring that only a specific number of skills be present across all team members combined. This approach has two consequences.
First, unless the minimum skill coverage is equal to the total number of skills, teams may differ in which skills they cover, which might be undesirable. For example, in a setting with four skills and a minimum skill coverage of two, one team might cover only the first two skills, while another covers only the last two. Although this solution is computationally valid, it may not align with the intended educational objectives. 
This issue can be mitigated by excluding non-essential skills and requiring all skills to be covered, while making sure that there are enough students to cover each skill marked as essential. Alternatively, future extensions of our model could prioritize certain skills to ensure that each team covers the most important ones. Second, our current model does not control the distribution of skills among team members.
In extreme cases, an all-rounder, i.e., a student who possesses all skills, could cover every skill for their team, permitting the other team members to cover no skill at all. 
While this may not be a problem in educational settings that emphasize individual development, it is less suitable for settings in which team members are expected to fill designated roles that require different skills. 
In our model, this can be controlled by reducing the likelihood of all-rounders being present during data collection, for example by raising the threshold for a student to be considered skilled or by capping the number of skills each student can cover.
Future work could also expand our model to include a more nuanced concept of skill coverage, such as minimizing skill overlap (i.e., avoiding assigning students with the same skills to the same team) or requiring explicit role assignments based on predefined team roles. 
In summary, future research could explore different ways to model skill coverage and role composition in a team, and evaluate their impact on team performance.

\section{Conclusion}
\label{sec:conclusion}

Algorithmic team formation is an area of growing interest in the education domain, where student teams are typically formed by human decision-makers, either teachers or students, often with suboptimal results. 
In this paper, we introduce the EDUCATIONAL TEAM FORMATION (EDU-TF) problem, which captures both top-down teacher requirements (such as skill coverage) and bottom-up student preferences for interpersonal compatibility within the teams.
We propose a modular and hierarchical optimization approach, which allows prioritizing the team formation objectives related to interpersonal preferences according to the needs of the particular educational setting, while ensuring skill coverage.
We evaluate four strategies (comprising ten sub-strategies) derived from this approach across eight real-world course datasets. Our results show that our method outperforms traditional, teacher-assigned teams, while maintaining computational efficiency.
Our findings also emphasize the importance of managing trade-offs, for example between maximizing preferred collaborations and minimizing undesired ones, highlighting the value of adaptable, context-aware algorithms in real classroom environments. 
This work advances the field of algorithmic team formation by bridging gaps between theoretical optimization and practical educational needs.
Future directions include scaling to larger datasets, exploring alternative skill coverage models, and experimenting with dynamic regrouping based on adaptations of the original solution to support evolving classroom dynamics.

\bibliographystyle{unsrtnat}
\bibliography{references}  






\end{document}